\documentclass[11pt]{article}
\usepackage{amsfonts}

\usepackage{makeidx}
\usepackage{amssymb}
\usepackage{graphicx}
\usepackage{amsmath}
\usepackage[T1]{fontenc}

0cm \voffset -1.5cm \DeclareMathOperator\Tr{Tr}
\newtheorem{theorem}{Theorem}

\newtheorem{corollary}[theorem]{Corollary}

\newtheorem{adefinition}[theorem]{Definition}

\newtheorem{aexample}[theorem]{Example}

\newtheorem{lemma}[theorem]{Lemma}

\newtheorem{proposition}[theorem]{Proposition}
\newtheorem{aremark}[theorem]{Remark}
\newenvironment{remark}{\begin{aremark}\rm}{\end{aremark}}

\numberwithin{equation}{section} \numberwithin{theorem}{section}
\newenvironment{proof}[1][Proof]{\textbf{#1.} }{\ \rule{0.5em}{0.5em}}
\makeatother

\begin{document}

\title{Non-Gaussian Limiting Laws for the Entries of Regular Functions of
the Wigner Matrices }
\author{A. Lytova, L. Pastur \\
Mathematical Division\\
Institute for Low Temperatures\\
Kharkiv, Ukraine}
\date{}
\maketitle

\begin{abstract}
This paper is a continuation of our paper \cite{L-Pa:09} in which
we proved the Central Limit Theorem for the matrix elements of
differential functions of the real symmetric random Gaussian
matrices (GOE). Here we consider the real symmetric random Wigner
matrices having independent (modulo symmetry conditions) but not
necessarily Gaussian entries. We show that in this case the matrix
elements of sufficiently smooth functions of these random matrices
have in general another limiting law which coincides essentially
with the probability law of matrix entries.
\end{abstract}

\section{Introduction}

We are interested in asymptotic properties of matrix elements $\varphi
_{jk}(M)$, $j,k=1,..,n$, $n\rightarrow \infty $, where $\varphi $ is a
smooth enough test-function, and $M$ is the Wigner matrix. We define the
Wigner real symmetric matrix as follows:
\begin{equation}
M_{n}=n^{-1/2}W_{n},\quad W_{n}=\{W_{jk}^{(n)}\in \mathbb{R}%
,\;W_{jk}^{(n)}=W_{kj}^{(n)}\}_{j,k=1}^{n},  \label{MW}
\end{equation}%
where $\{W_{jk}^{(n)}\}_{1\leq j\leq k\leq n}$ are independent random
variables satisfying
\begin{equation}
\mathbf{E}\{W_{jk}^{(n)}\}=0,\quad \mathbf{E}\{(W_{jk}^{(n)})^{2}\}=w^{2}(1+%
\delta _{jk}).  \label{Wmom12}
\end{equation}%
The case of the Gaussian random variables
obeying (\ref{Wmom12}) corresponds to the Gaussian Orthogonal Ensemble (GOE)
(see e.g. \cite{Me:91}):
\begin{equation}
\widehat{M_{n}}=n^{-1/2}\widehat{W}_{n},\quad \widehat{W}_{n}=\{\;\widehat{W}%
_{jk}=\widehat{W}_{kj}\in \mathbb{R},\;\widehat{W}_{jk}\in \mathcal{N}%
(0,w^{2}(1+\delta _{jk}))\}_{j,k=1}^{n}.  \label{GOE}
\end{equation}%
We will assume in what follows additional conditions on distributions of $%
W_{jk}^{(n)}$, mostly in the form of existence of certain moments of $%
W_{jk}^{(n)}$, whose order will depend on the problem under study.

In our paper \cite{L-Pa:09} we have considered matrix elements of functions
of the GOE matrices and have proved the following facts.

\begin{theorem}
\label{t:covGE} Let $\widehat{M_{n}}$ be the GOE matrix (\ref{GOE}), and $%
\varphi _{1,2}:\mathbb{R\rightarrow C}$ be bounded test functions with
bounded derivative. Denote
\begin{equation}
({\varphi _{i}(\widehat{M_{n}})})_{jj}^{\circ }=(\varphi _{i}(\widehat{M_{n}}%
))_{jj}-\mathbf{E}\{(\varphi _{i}(\widehat{M_{n}}))_{jj}\}  \label{fjjc}
\end{equation}%
and
\begin{equation*}
\mathbf{Cov}\{(\varphi _{1}(\widehat{M_{n}}))_{jj},(\varphi _{2}(\widehat{%
M_{n}}))_{jj}\}=\mathbf{E}\{(\varphi _{1}(\widehat{M_{n}}))_{jj}({\varphi
_{2}(\widehat{M_{n}})})_{jj}^{\circ }\}.
\end{equation*}%
Then we have for any $j_{n}=1,...,n$
\begin{align}
\lim_{n\rightarrow \infty }& n\,\mathbf{Cov}\Big\{ \big( \varphi _{1}(%
\widehat{M_{n}})\big) _{j_{n}j_{n}},\big( \varphi _{2}(\widehat{M_{n}}%
)\big )_{j_{n}j_{n}}\Big\}  \notag \\
& =\int_{-2w}^{2w}\int_{-2w}^{2w}\Delta \varphi _{1}\Delta
\varphi _{2}\rho _{sc}(\lambda _{1})\rho _{sc}(\lambda _{2})d\lambda
_{1}d\lambda _{2},  \label{Cjjphi}
\end{align}%
where
\begin{equation}
\triangle \varphi =\varphi (\lambda _{1})-\varphi (\lambda _{2}),
\label{deltaf}
\end{equation}%
$\rho _{sc}$ is the density of the semicircle law%
\begin{equation}
\rho _{sc}(\lambda )=(2\pi w^{2})^{-1}(4w^{2}-\lambda ^{2})_{+}^{1/2},
\label{rhosc}
\end{equation}%
and $x_{+}=\max \{0,x\}$.
\end{theorem}

\begin{theorem}
\label{t:cltGE} Let $\widehat{M_{n}}$ be the GOE matrix (\ref{GOE}), and $%
\varphi :\mathbb{R\rightarrow R}$ be bounded function with bounded
derivative. Then for any $j_{n}=1,...,n$ the random variable $\sqrt{n}%
\varphi ^{\circ }(\widehat{M_{n}})_{j_{n}j_{n}}$ converges in distribution
to the Gaussian random variable with zero mean and the variance
\begin{equation}
V_{d}^{GOE}[\varphi ]=\int_{-2w}^{2w}\int_{-2w}^{2w}(\Delta
\varphi )^{2}\rho _{sc}(\lambda _{1})\rho _{sc}(\lambda _{2})d\lambda
_{1}d\lambda _{2}.  \label{Varjj}
\end{equation}
\end{theorem}

In the present paper we prove a counterpart of Theorems \ref{t:covGE} and %
\ref{t:cltGE} for the Wigner matrices. In particular, we show in Theorem \ref%
{t:Cov} below that in this case the r.h.s. of (\ref{Cjjphi}) has an
additional term proportional to the fourth cumulant $\kappa _{4}$ of
non-diagonal entries (see (\ref{cums}) for the definition). This result is
in accordance with that for centered linear eigenvalue statistics%
\begin{equation}
\mathcal{N}_{n}^{\circ }[\varphi ]=\Tr\varphi^{\circ } (M_{n})=\sum_{j=1}^{n}\varphi_{jj}
^{\circ }(M_{n})  \label{ls}
\end{equation}%
of the Wigner matrices (see Theorem 3.6 of \cite{Ly-Pa:08}). On the other
hand the individual matrix elements $\sqrt{n}\varphi_{jj} ^{\circ }(M)$ do
not satisfy in general the Central Limit Theorem (CLT). In Theorem \ref%
{t:clt} of this paper we find limiting probability law for $\sqrt{n}\varphi_{jj}
^{\circ }(M)$ which is not Gaussian in general but rather that of
the sum of the probability law of entries of $M_{n}$ modulo a certain
rescaling. To obtain the CLT, one has to impose an integral condition on the
test function, i.e., the set of test functions for which we have the CLT has
the codimension one.

Our results of \cite{L-Pa:09} and of this paper can be viewed as analogs of
the E. Borel theorem on the limiting probability law of entries of orthogonal
matrices of size $n$ as $n \to \infty$ (see e.g. \cite{Di-Co:03}).

\textit{Convention}: We will use letter $C$ for an absolute constant that
does not depend on $j$, $k$, and $n$, and may be distinct on different
occasions.

\section{Technical Means}

To make the paper self-consistent, we present here several technical facts
that will be often used below. For the proof of these facts see e.g. \cite%
{KKP:96,Ly-Pa:08}

We start from the generalized Fourier transform, in fact the $\pi /2$
rotated Laplace transform (see e.g. \cite{Ti:86}, Sections 1.8-9 for its
definition).

\begin{proposition}
\label{p:Four} Let $f:\mathbb{R}_{+}\rightarrow \mathbb{C}$ be a locally
Lipshitzian and such that for some $\delta >0$%
\begin{equation}
\sup_{t\geq 0}e^{-\delta t}|f(t)|<\infty ,  \label{supde}
\end{equation}%
and let $\widetilde{f}:\{z\in \mathbb{C}:\Im z<-\delta \}\rightarrow \mathbb{%
C}$ be its \emph{generalized Fourier transform}
\begin{equation}
\widetilde{f}(z)=i^{-1}\int_{0}^{\infty }e^{-izt}f(t)dt.  \label{Fur}
\end{equation}%
The inversion formula is given by
\begin{equation}
f(t)=\frac{i}{2\pi }\int_{{L}}e^{izt}\widetilde{f}(z)dz,\;t\geq 0,
\label{Furinv}
\end{equation}%
where ${L}=(-\infty -i\varepsilon ,\infty -i\varepsilon )$, $\varepsilon
>\delta ,$ and the principal value of the integral at infinity is used.

Denote for the moment the correspondence between functions and their
generalized Fourier transforms as $f\leftrightarrow \widetilde{f}$. Then we
have:

\begin{enumerate}
\item[(i)] $\quad \int_{0}^{t}f(\tau )d\tau \leftrightarrow (iz)^{-1}%
\widetilde{f}(z);$

\item[(ii)] $\quad \int_{0}^{t}f_{1}(t-\tau )f_{2}(\tau )d\tau :=(f_{1}\ast
f_{2})(t)\leftrightarrow i \widetilde{f_{1}}(z)\widetilde{f_{2}}(z);$

\item[(iii)] if $P$, $Q$, and $R$ are differentiable, and $R(0)=0$, then the
equation
\begin{equation}
P(t)+\int_{0}^{t}dt_{1}%
\int_{0}^{t_{1}}Q_{{}}(t_{1}-t_{2})P(t_{2})dt_{2}=R(t),\;t\geq 0,
\label{intrel}
\end{equation}%
has a unique differentiable solution
\begin{equation}
P(t)=-\int_{0}^{t}T_{{}}(t-t_{1})R^{\prime }(t_{1})dt_{1},  \label{solut}
\end{equation}%
where
\begin{equation}
T\leftrightarrow (z+\widetilde{Q})^{-1}  \label{TtSt}
\end{equation}%
provided by
\begin{equation}
z+\widetilde{Q}(z)\neq 0,\;\Im z<0.  \label{condQ1}
\end{equation}%
%
%
%
%
%
%
%
\end{enumerate}
\end{proposition}

\noindent The next proposition presents simple facts of linear algebra

\begin{proposition}
\label{p:Duh} Let $M$ and $M^{\prime }$ be $n\times n$ matrices and $t\in
\mathbb{R}$. Then we have the following:

\begin{enumerate}
\item[(i)] the Duhamel formula
\begin{equation}
e^{(M+M^{\prime })t}=e^{Mt}+\int_{0}^{t}e^{M(t-s)}M^{\prime }e^{(M+M^{\prime
})s}ds.  \label{Duh}
\end{equation}

\item[(ii)] if for a real symmetric $M$
\begin{equation}
U(t)=e^{itM},\;t\in \mathbb{R},  \label{U}
\end{equation}%
then $U(t)$ is a symmetric unitary matrix and
\begin{equation}
U(t_{1})U(t_{2})=U(t_{1}+t_{2}),\quad ||U(t)||=1,\quad
\sum_{j=1}^{n}|U_{jk}(t)|^{2}=1,  \label{norU}
\end{equation}
so that
\begin{align}
&|U_{jk}(t)|\leq 1,\quad \sum_{k=1}^n|U_{jk}(t)|\leq n^{1/2},\quad
\sum_{k=1}^{n}|U_{j_{}k}(t_{1})U_{j^{\prime }k}(t_{2})|\leq 1;  \label{SUjk<}
\end{align}

\item[(iii)]
\begin{equation}
D_{jk}U_{ab}(t)=i\beta _{jk}\left (U_{aj}\ast U_{bk}+U_{bj}\ast
U_{ak}\right)(t),  \label{ParU}
\end{equation}%
where
\begin{align}
&D_{jk}=\partial /\partial M_{jk},  \label{Djk} \\
&\beta _{jk}=(1+\delta _{jk})^{-1}=1-\delta_{jk}/2,  \label{beta}
\end{align}%
the symbol "$\ast $" is defined in Proposition \ref{p:Four} (ii), and
\begin{equation}
|D_{jk}^{l}U_{ab}(t)|\leq c_{l}|t|^{l},\quad c_l=2^l/l!.  \label{ocdlu}
\end{equation}
\end{enumerate}
\end{proposition}

Now a generalization of a property of the Gaussian random variable $\xi$ of
zero mean and of variance $w^2$ according to which if $\Phi :\mathbb{R}
\rightarrow \mathbb{C}$ is a differentiable function with polynomially
bounded derivative, then
\begin{equation}
\mathbf{E}\{\xi \, \Phi(\xi )\}=\mathbf{E}\{\xi^{2}\}\mathbf{E}%
\{\Phi^{\prime}\}.  \label{diffga}
\end{equation}%
Formula (\ref{diffga}) is a particular case of more general formula. To
write it we recall some definitions. If a random variable $\xi $ has a
finite $p$th absolute moment, $p\geq 1$, then we have the expansions
\begin{equation*}
f(t):=\mathbf{E}\{e^{it\xi }\}=\sum_{j=0}^{p}\frac{\mu _{j}}{j!}%
(it)^{j}+o(t^{p}),
\end{equation*}%
and
\begin{equation}
l(t):=\log \mathbf{E}\{e^{it\xi }\}=\sum_{j=0}^{p}\frac{\kappa _{j}}{j!}%
(it)^{j}+o(t^{p}),\quad t\rightarrow 0,  \label{lt}
\end{equation}%
where $"\log "$ denotes the principal branch of logarithm, the coefficients
in the expansion of $f$ are the moments $\{\mu _{j}\}$ of $\xi $, and the
coefficients in the expansion of $l$ are the cumulants $\{\kappa _{j}\}$ of $%
\xi $. For small $j$ one easily expresses $\kappa _{j}$ via $\mu _{1},\mu
_{2},\dots ,\mu _{j}$. In particular,
\begin{align}
\kappa _{1}& =\mu _{1},\quad \kappa _{2}=\mu _{2}-\mu _{1}^{2}=\mathbf{Var}%
\{\xi \},\quad \kappa _{3}=\mu _{3}-3\mu _{2}\mu _{1}+2\mu _{1}^{3},
\label{cums} \\
\kappa _{4}& =\mu _{4}-3\mu _{2}^{2}-4\mu _{3}\mu _{1}+12\mu _{2}\mu
_{1}^{2}-6\mu _{1}^{4},...  \notag
\end{align}%
In general%
\begin{equation}
\kappa _{j}=\sum_{\lambda }c_{\lambda }\mu _{\lambda },  \label{kapmu}
\end{equation}%
where the sum is over all additive partitions $\lambda $ of the set $%
\{1,\dots ,j\}$, $c_{\lambda }$ are known coefficients and $\mu _{\lambda
}=\prod_{l\in \lambda }\mu _{l}$, see e.g. \cite{Pr-Ro:69}. We have

\begin{proposition}
\label{l:difgen} Let $\xi $ be a random variable such that $\mathbf{E}\{|\xi
|^{p+2}\}<\infty $ for a certain non-negative integer $p$. Then for any
function $\Phi :\mathbb{R}\rightarrow \mathbb{C}$ \ of the class $C^{p+1}$
with bounded partial derivatives $\Phi ^{(l)}$, $l=1,..,p+1$, we have
\begin{equation}
\mathbf{E}\{\xi \Phi (\xi )\}=\sum_{l=0}^{p}\frac{\kappa _{l+1}}{l!}\mathbf{E%
}\{\Phi ^{(l)}(\xi )\}+\varepsilon _{p},  \label{difgen}
\end{equation}%
where
\begin{equation}
|\varepsilon _{p}|\leq C_{p}\mathbf{E}\{|\xi |^{p+2}\}\sup_{t\in \mathbb{R}%
}|\Phi ^{(p+1)}(t)|,\,\,\,C_{p}\leq \frac{1+(3+2p)^{p+2}}{(p+1)!}.
\label{b3}
\end{equation}
If the characteristic function $\mathbf{E}\{e^{it|\xi| }\}$ is whole, and $%
\Phi \in C^\infty$, then
\begin{equation}
\mathbf{E}\{\xi \Phi (\xi )\}=\sum_{l=0}^{\infty}\frac{\kappa _{l+1}}{l!}%
\mathbf{E}\{\Phi ^{(l)}(\xi )\}  \label{difinf}
\end{equation}
provided that for some $a>0$
\begin{equation}
|\mathbf{E}\{\Phi ^{(l)}(\xi )\}|\leq a^l,  \label{al}
\end{equation}
and for some $R=ca,$ $c>1$, $\kappa_l$ satisfy the condition:
\begin{equation}
\sum_{l=0}^{\infty}\frac{|\kappa _{l+1}|R^l}{l!}<\infty.  \label{kap<}
\end{equation}
\end{proposition}

\section{Covariance of Matrix Elements}

We show first that if $M$ is the Wigner matrix with uniformly bounded sixth
moments of its entries, and the test-function $\varphi$ is essentially of
class $\mathbf{C}^3$, then the variance of $\sqrt{n}\varphi _{jj}^{ }(M)$ is
of the order $O(1)$ as $n\rightarrow\infty$. We have

\begin{lemma}
\label{l:Var} Let $M=n^{-1/2}W$ be the real symmetric Wigner matrix (\ref{MW}%
) -- (\ref{Wmom12}). Assume that:

(i) the third moments of its entries do not depend on $j$, $k$, and $n$:%
\begin{equation}  \label{mu3}
\mu _{3}=\mathbf{E}\big\{ ( W_{jk}^{(n)}){}^{3}\big\};
\end{equation}

(ii) the sixth moments are uniformly bounded:%
\begin{equation}  \label{w6<}
w_{6}:=\sup_{n\in \mathbb{N}}\max_{1\leq j,k\leq n}\mathbf{E}\big\{ (
W_{jk}^{(n)}){}^{6}\big\} <\infty.
\end{equation}

\noindent Then for any test-function $\varphi :\mathbb{R\rightarrow C}$,
whose Fourier transform
\begin{equation}
\widehat{\varphi }(t)=\frac{1}{2\pi }\int e^{-it\lambda }\varphi (\lambda
)\;d\lambda  \label{FT}
\end{equation}%
satisfies the condition
\begin{equation}
\int(1+|t|)^{3}|\widehat{\varphi }(t)|dt<\infty,  \label{F3<}
\end{equation}%
we have the bound
\begin{align}
\mathbf{Var}\{\sqrt{n}\varphi _{jj}(M)\}:&=\mathbf{E}\{|\sqrt{n}%
\varphi^\circ_{jj}(M)|^2\}  \notag \\
&\leq C\Big( \int (1+|t|)^{3}|\widehat{\varphi }(t)|dt\Big) ^{2}.
\label{VarF}
\end{align}
\end{lemma}

\begin{proof}
It follows from the Fourier inversion formula
\begin{equation}
\varphi (\lambda )=\int e^{i\lambda t}\widehat{\varphi }(t)dt  \label{Finv}
\end{equation}%
and the spectral theorem for symmetric matrices that
\begin{equation}
\varphi _{jj}^{\circ }(M)=\int \widehat{\varphi }(t)U_{jj}^{\circ }(t)dt,
\label{phijj}
\end{equation}%
where $U$ is defined in (\ref{U}). This and the Schwarz inequality yield
\begin{align}
\mathbf{Var}\{\varphi _{jj}(M)\}&=\int\int \mathbf{E}\{U_{jj}(t_1)U_{jj}^{%
\circ }(t_2)\}\widehat{\varphi }(t_1)\widehat{\varphi }(t_2)d t_{1}dt_{2}
\label{Varphi} \\
&\leq\bigg(\int \mathbf{Var}^{1/2}\{U_{jj}(t)\}|\widehat{\varphi }(t)|dt%
\bigg)^2.  \notag
\end{align}
Now (\ref{VarF}) follows from the estimate
\begin{equation}
\mathbf{Var}\{U_{jj}(t)\} \leq C(1+|t|)^6/n  \label{WUjj<}
\end{equation}%
proved in Lemma \ref{l:main} below (see Appendix, (\ref{Ujjv})), and
condition (\ref{F3<}).
\end{proof}

\begin{theorem}
\label{t:Cov} Let $M=n^{-1/2}W$ be the real symmetric Wigner matrix (\ref{MW}%
) -- (\ref{Wmom12}). Assume that the third and fourth moments do not depend
on $j$, $k$, and $n$:%
\begin{equation}  \label{mu4}
\mu _{3}=\mathbf{E}\big\{ ( W_{jk}^{(n)}){}^{3}\big\} ,\quad\mu _{4}=\mathbf{%
E}\big\{ ( W_{jk}^{(n)}){}^{4}\big\} ,
\end{equation}
and the sixth absolute moments are uniformly bounded (see (\ref{w6<})).
\noindent Let $\varphi_{1,2} :\mathbb{R\rightarrow C}$ be the
test-functions, whose Fourier transforms $\widehat{\varphi}_{1,2}$ (\ref{FT}%
) satisfy (\ref{F3<}). Then we have for any $j=j_{n}\in \lbrack 1,n]$:
\begin{align}
\lim_{n\rightarrow\infty}n\mathbf{Cov}\{(\varphi_1 (M))_{jj}, (\varphi_2
(M))_{jj}\} =&\int_{-2w}^{2w}\int_{-2w}^{2w}\triangle\varphi_1{%
\triangle\varphi_2} \rho _{sc}(\lambda _{1})\rho _{sc}(\lambda _{2})d\lambda
_{1}d\lambda _{2}  \notag  \label{Coff} \\
&+\frac{\kappa_4}{w^{8}}\prod_{i=1}^2\int_{-2w}^{2w}\varphi_i(\lambda
)(w^2-\lambda^2)\rho _{sc}(\lambda)d\lambda,
\end{align}
where $\triangle\varphi$ is defined in (\ref{deltaf}), and
\begin{equation}
\kappa_4=\mu_4-3w^4  \label{k4}
\end{equation}
is the fourth cumulant of the off-diagonal entries (see (\ref{cums})).

In particular,
\begin{align}
V_d^W[\varphi]:&=\lim_{n\rightarrow\infty}\mathbf{Var}\{\sqrt{n}\varphi
_{jj}(M)\}  \notag \\
&=V_d^{GOE}[\varphi]+\frac{\kappa_4}{w^{8}}\Big| \int_{-2w}^{2w}\varphi(%
\mu)(w^2-\mu^2)\rho _{sc}(\mu)d\mu\Big|^2  \label{VWd}
\end{align}
with $V_d^{GOE}[\varphi]$ of (\ref{Varjj}).
\end{theorem}

\begin{remark}
\label{r:odd} \textit{(i).} If $\varphi$ is odd, then $V_d^{W}[%
\varphi]=V_d^{GOE}[\varphi]$.

\textit{(ii).} Note that we choose here the Wigner matrix so that its first
two moments matches the first two moments of the GOE matrix (see (\ref%
{Wmom12})). This fact allows to use known properties of GOE and lies at the
basis of interpolation procedure widely used in the proof of Lemma \ref%
{l:main} below. In fact this condition is pure technical one, and we can
replace condition (\ref{Wmom12}) with more general one and consider Wigner
matrix $\widetilde{M}=n^{-1/2}\widetilde{W}$, satisfying
\begin{align}
&\mathbf{E}\{\widetilde{W}^{(n)}_{jk}\}=0,\quad 1\leq j\leq k\leq n,
\label{mom12} \\
&\mathbf{E}\{(\widetilde{W}^{(n)}_{jk})^{2}\}=w^{2},\; j\neq k,\quad \mathbf{%
E}\{(\widetilde{W}^{(n)}_{jj})^{2}\}=w_{2}w^{2},\; w_2>0.  \notag
\end{align}%
In this case there arise additional terms in (\ref{Coff}) and (\ref{VWd})
proportional to $w_2-2$. In particular, we have for the corresponding
limiting variance
\begin{align}
V_d^{\widetilde{W}}[\varphi]=V_d^W[\varphi]+(w_2-2)w^{-2}\bigg| %
\int_{-2w}^{2w}\varphi(\mu)\mu\rho _{sc}(\mu)d\mu\bigg|^2,  \label{VWd2}
\end{align}
where $V_d^W[\varphi]$ is given by (\ref{VWd}).
\end{remark}

\begin{proof}
Let us write the covariation in the form (cf (\ref{Varphi}))
\begin{equation}
n\mathbf{Cov}\{(\varphi_1 (M))_{jj}, (\varphi_2 (M))_{jj}\}= \int\int
C_n(t_1,t_2)\widehat{\varphi}_1(t_1)\widehat{\varphi}_2(t_2)dt_1dt_2,
\label{nCov}
\end{equation}
where
\begin{equation*}
C_n(t_1,t_2)=n\mathbf{E}\{U_{jj}(t_1)U_{jj}^{\circ}(t_2)\},
\end{equation*}
and $U$ is defined in (\ref{U}). It follows from the Schwarz inequality and (%
\ref{WUjj<}) that
\begin{equation*}
C_n(t_1,t_2)\leq C(1+|t_1|)^3(1+|t_2|)^3.
\end{equation*}
This, (\ref{F3<}), and (\ref{nCov}) imply that it suffices to show that
there are converging subsequences $\{C_{n_i}\}$ and function $Cov$ such that
we have for any converging subsequence $\{C_{n_i}\}$ and any $T>0$
\begin{equation}
\lim_{i\rightarrow\infty}C_{n_i}(t_1,t_2)=Cov(t_1,t_2)  \label{Cov}
\end{equation}
uniformly on the square $S_T=\{(t_1,t_2)\in\mathbb{R}:|t_1|\leq T,|t_2|\leq
T\}$, and that plugging $Cov$ in the r.h.s. of (\ref{nCov}) we get the
r.h.s. of (\ref{Coff}).

Show first that the derivatives $\partial C_n/\partial t_{i}$, $i=1,2$, are
bounded on $S_T$ uniformly in $n$. We have
\begin{align}
\frac{\partial C_n}{\partial t_{1}}=in\mathbf{E}\{(MU)_{jj}(t_1)U_{jj}^{%
\circ}(t_2)\}=i\sqrt{n}\sum_{k=1}^n \mathbf{E}\{{W}^{(n)}_{jk}%
\Phi_{jk}(t_1,t_2)\},  \label{dtC}
\end{align}
where
\begin{equation}
\Phi_{jk}(t_1,t_2)=U_{jk}(t_1)U_{jj}^{\circ}(t_2).  \label{Fjk}
\end{equation}
A simple algebra based on (\ref{ocdlu}) allows to obtain
\begin{equation}
|D^{l}_{jk}\Phi_{jk}(t_3,t_2)|=O(1),\quad n\rightarrow\infty,  \label{Fjk<}
\end{equation}
uniformly in $(t_1,t_2)\in S_T$. Now applying differentiation formula (\ref%
{difgen}) with $\Phi=\Phi_{jk}$ and $p=2$ to every term of the r.h.s. of (%
\ref{dtC}), we get
\begin{align}
\frac{\partial C_n}{\partial t_{1}}=iw^{2}\sum_{k=1}^n \beta_{jk}^{-1}%
\mathbf{E}\{D_{jk}\Phi_{jk}(t_1,t_2)\}+\frac{i\mu_3}{\sqrt{n}}\sum_{k=1}^n
\mathbf{E}\{D^2_{jk}\Phi_{jk}(t_1,t_2)\}+\varepsilon_{2}(t_{1},t_{2}),
\label{dtC2}
\end{align}
where in view of (\ref{b3}) and (\ref{Fjk<})
\begin{equation}
|\varepsilon _{2}(t_1,t_2)|\leq \frac{C_{2}\mu_4}{n}\sum_{k=1}^{n}\sup_{M\in
\mathcal{S}_{n}}\big|D_{jk}^{3}\Phi _{jk}(t_1,t_2)\big|=O(1),\quad
n\rightarrow\infty.  \label{e2}
\end{equation}
Here $\mathcal{S}_{n}$ is the set of $n\times n$ real symmetric matrices.
Now it follows from (\ref{ParU}) and (\ref{Fjk}) that every term of $%
D^2_{jk}\Phi_{jk}$ contains $U_{jk}$ (see (\ref{T2=})). Taking into account
that $\sum_{k=1}^n|U_{jk}|\leq n^{1/2}$ (see (\ref{SUjk<})), we see that the
second term on the r.h.s. of (\ref{dtC2}) is of the order $O(1)$, $%
n\rightarrow\infty$ uniformly in $(t_1,t_2)\in S_T$.  At last, using (\ref%
{ParU}) we get for the first term on the r.h.s. of (\ref{dtC2}):
\begin{align*}
T_1^{(n)}(t_1,t_2):&=w^{2}\sum_{k=1}^n \beta_{jk}^{-1}\mathbf{E}%
\{D_{jk}\Phi_{jk}(t_1,t_2)\} \\
&=iw^2\sum_{k=1}^{n} \mathbf{E}\big\{%
(U_{jj}*U_{kk}+U_{jk}*U_{jk})(t_1)U_{jj}^{\circ}(t_2)
+2U_{jk}(t_1)(U_{jj}*U_{jk})(t_2)\big\} \\
&=iw^2\mathbf{E}\big\{\big[(U_{jj}*nv_n)(t_1)+t_1U_{jj}(t_1)\big]%
U_{jj}^{\circ}(t_2)+2\int_0^{t_2}U_{jj}(t_1+t_4)U_{jj}(t_2-t_4)dt_4\big\},
\end{align*}
where we denote
\begin{equation}  \label{vn=}
v_n(t):=n^{-1}\Tr U(t)=n^{-1}\sum_{k=1}^n U_{kk}(t),\quad |v_n(t)|\leq 1.
\end{equation}
Since
\begin{equation*}
n\mathbf{E}\{v_nU_{jj}U_{jj}^{\circ}\}=\mathbf{E}\{v_n\}C_n+n\mathbf{E}%
\{v_n^{\circ}U_{jj}U_{jj}^{\circ}\}
\end{equation*}
and by (\ref{Ujjv}) -- (\ref{Varu<}) of Lemma \ref{l:main} below
\begin{equation}  \label{VarUv}
\mathbf{Var}\{U_{jj}\}=O(n^{-1}), \quad\mathbf{Var}\{v_n\}=O(n^{-2}),\quad
n\rightarrow\infty,
\end{equation}
uniformly in $(t_1,t_2)\in S_T$, then we finally have
\begin{align}
T_1^{(n)}(t_1,t_2)=&iw^2\int_0^{t_1}\mathbf{E}\{v_n(t_4)%
\}C_n(t_1-t_4,t_2)dt_4  \notag \\
&+2iw^2\int_0^{t_2}\mathbf{E}\{U_{jj}(t_1+t_4)\}\mathbf{E}%
\{U_{jj}(t_2-t_4)\}dt_4+r_n(t_1,t_2),  \label{T1=}
\end{align}
where
\begin{align}
r_n(t_1,t_2)=&iw^2\mathbf{E}\big\{\big[(U_{jj}*nv_n^%
\circ)(t_1)+t_1U_{jj}(t_1)\big] U_{jj}^{\circ}(t_2)\big\}  \notag \\
&+2\int_0^{t_2}\mathbf{E}\big\{U_{jj}(t_1+t_4)U_{jj}^\circ(t_2-t_4)\big\}%
dt_4=O(n^{-1/2}),\quad n\rightarrow\infty,  \label{r}
\end{align}%
and see that $T_1^{(n)}(t_1,t_2)=O(1)$,  $n\rightarrow\infty$ uniformly in $%
(t_1,t_2)\in S_T$.

It follows from the above that the derivatives $\partial C_n/\partial t_{i}$%
, $i=1,2$, are  bounded on $S_T$ uniformly in $n$. Hence, there are
converging subsequences $\{C_{n_i}\}$ and function $Cov$ (depending on
subsequence) such that (\ref{Cov}) holds. Now we derive an integral equation
for $Cov$ showing that $Cov$ is the same for every converging subsequences $%
\{C_{n_i}\}$ and leading via (\ref{nCov}) to (\ref{Coff}).

It follows from the Duhamel formula (\ref{Duh}) that
\begin{align*}
C_n(t_1,t_2)&=\int_0^{t_1}in^{1/2}\sum_{k=1}^n \mathbf{E}\{{W}%
^{(n)}_{jk}\Phi_{jk}(t_3,t_2)\}dt_3
\end{align*}
with $\Phi_{jk}(t_3,t_2)$ given by (\ref{Fjk}). We see that the integrand
here coincides with the r.h.s. of (\ref{dtC}). Hence, applying
differentiation formula (\ref{difgen}) with $p=3$, we get (cf (\ref{dtC2})
-- (\ref{e2})):
\begin{align}
C_n(t_1,t_2)=\int_0^{t_1}i\bigg[\sum_{l=1}^3
T_l^{(n)}(t_3,t_2)+\varepsilon_3(t_3,t_2)\bigg]dt_3,  \label{Covn}
\end{align}
where
\begin{align}
&T_l^{(n)}(t_3,t_2)=\frac{1}{l!n^{(l-1)/2}}\sum_{k=1}^{n}\kappa _{l+1,jk}%
\mathbf{E}\big\{D_{jk}^{l}\Phi_{jk}(t_3,t_2)\big\},\quad l=1,2,3,
\label{Tl1}
\end{align}
$\kappa _{l,jk}$ is the $l$th cumulant of ${W}^{(n)}_{jk}$:
\begin{equation}  \label{kap1234}
\kappa _{1,jk}=0,\quad\kappa _{2,jk}=w^2\beta^{-1}_{jk},\quad \kappa
_{3,jk}=\mu_3,\quad \kappa _{4,jk}=\kappa _{4}-9\delta_{jk}w^4,
\end{equation}
(see (\ref{Wmom12}), (\ref{cums}), and (\ref{k4})), and in view of (\ref{b3}%
) and (\ref{Fjk<})
\begin{equation}
|\varepsilon _{3}(t_3,t_2)|\leq \frac{C_{3}w^{5/6}_{6}}{n^{3/2}}%
\sum_{k=1}^{n}\sup_{M\in \mathcal{S}_{n}}\big|D_{jk}^{4}\Phi _{jk}(t_3,t_2)%
\big|=O( n^{-1/2}),\quad n\rightarrow\infty.  \label{e3clt4}
\end{equation}
We see that $T^{(n)}_1$ of (\ref{Tl1}) is given by (\ref{T1=}) -- (\ref{r}),
where by  (\ref{vt}) -- (\ref{Ujjv}) of Lemma \ref{l:main} below
\begin{equation}
\lim_{n\rightarrow\infty}\mathbf{E}\{v_n(t)\}=\lim_{n\rightarrow\infty}%
\mathbf{E}\{U_{jj}(t)\}=v(t)=\int e^{it\lambda}\rho_{sc}(\lambda)d\lambda.
\label{limU}
\end{equation}
This and  (\ref{Cov}) yield
\begin{align}
\lim_{i\rightarrow\infty}T_1^{(n_i)}(t_3,t_2)=iw^2%
\int_0^{t_3}v(t_4)Cov(t_3-t_4,t_2)dt_4 +2iw^2\Phi(t_3,t_2),  \label{T1}
\end{align}
where
\begin{align}
\Phi(t_3,t_2)=i^{-1}\int_{-2w}^{2w}\int_{-2w}^{2w}e^{it_3\lambda}\frac{%
e^{it_2\lambda}-e^{it_2\mu}}{\lambda-\mu}\rho _{sc}(\lambda)\rho
_{sc}(\mu)d\lambda d\mu.  \label{Phi}
\end{align}
Consider now $T_2^{(n)}$ of (\ref{Tl1}), and show that
\begin{align}
\lim_{n\rightarrow\infty}T_2^{(n)}(t_3,t_2)=0.  \label{T2}
\end{align}
We have by (\ref{ParU}) and (\ref{Tl1}) with $l=2$:
\begin{align}
T_2^{(n)}(t_3,t_2)=-\frac{\mu_3}{n^{1/2}}\sum_{k=1}^n \beta_{jk}^2 \mathbf{E}%
\{& (U_{jk}* U_{jk}* U_{jk}+3U_{jj}* U_{jk}* U_{kk})(t_3)U_{jj}^\circ(t_{2})
\notag \\
&+2(U_{jk}* U_{jk}+ U_{jj}*U_{kk})(t_3)(U_{jj}*U_{jk})(t_2)  \notag \\
&+U_{jk}(t_3)(3U_{jj}*U_{jk}*U_{jk}+U_{jj}*U_{jj}*U_{kk})(t_2)\}.
\label{T2=}
\end{align}
It follows from (\ref{SUjk<}) that the contribution of the terms containing $%
U_{jk}U_{jk}U_{jk}$  is of the order $O(n^{-1/2})$, $n\rightarrow\infty$.
Besides,  since $n^{-1/2}\sum_{k=1}^{n}|U_{kk}(t)U_{jk}(\tau)|\leq 1$, then
by the Schwarz inequality and (\ref{WUjj<}) the contribution  of the terms
containing $U_{jj}^\circ(t_{2})$  is also of the order $O(n^{-1/2})$, $%
n\rightarrow\infty$. So  we are left with
\begin{equation*}
-\frac{\mu_3}{n^{1/2}}\sum_{k=1}^n \mathbf{E}%
\{2(U_{jj}*U_{kk})(t_3)(U_{jj}*U_{jk})(t_2)
+U_{jk}(t_3)(U_{jj}*U_{jj}*U_{kk})(t_2)\}.
\end{equation*}%
Here by (\ref{v1}) of Lemma \ref{l:main} below
\begin{align}
&\mathbf{Var}\Big\{{n^{-1/2}}\sum_{k=1}^n U_{jk}(\tau_1)U_{kk}(\tau_2)\Big\}%
=O(n^{-1/2}),\;n\rightarrow\infty,  \label{jkkk1} \\
&\lim_{n\rightarrow\infty}\mathbf{E}\Big\{{n^{-1/2}}\sum_{k=1}^n
U_{jk}(\tau_1)U_{kk}(\tau_2)\Big\}=0,  \label{jkkk2}
\end{align}%
so that
\begin{align*}
\lim_{n\rightarrow\infty}\mathbf{E}\Big\{{n^{-1/2}}\sum_{k=1}^n
U_{jk}(\tau_1)U_{kk}(\tau_2)U_{jj}(\tau_3)U_{jj}(\tau_4)\Big\}=0,
\end{align*}%
and we get (\ref{T2}).

Consider now $T_3^{(n)}$ of (\ref{Tl1}). We have
\begin{align*}
T_3^{(n)}(t_3,t_2)=\frac{\kappa_4}{6n}\sum_{k=1}^n\mathbf{E}\{
D_{jk}^{3}(U_{jk}(t_3)U_{jj}^\circ(t_2))\}+O(n^{-1}),\;{n\rightarrow\infty},
\end{align*}
where we replaced $\kappa_{4,jk}$ of (\ref{kap1234}) with $\kappa_{4}$ of (%
\ref{k4}) with the error term of the order $O(n^{-1})$, ${n\rightarrow\infty}
$. It follows now from (\ref{SUjk<}) -- (\ref{ParU}) that the contribution
to $T_3^{(n)}$ due to any term of
\begin{equation*}
n^{-1}\sum_{k=1}^{n}D_{jk}^{3}(U_{jk}(t_3)U_{jj}^\circ(t_2))
\end{equation*}%
containing at least one off-diagonal element $U_{jk}$ is of the order $%
O(n^{-1/2})$, ${n\rightarrow\infty}$. Besides, we have by (\ref{ocdlu}) and (%
\ref{WUjj<}) that the term $\mathbf{E}\{U_{jj}^\circ(t_2)n^{-1}%
\sum_{k=1}^{n}D_{jk}^{3}U_{jk}(t_3)\}$ is of the order $O(n^{-1/2})$, ${%
n\rightarrow\infty}$, too. Thus, we are left with terms, containing only
diagonal non-centered elements of $U$. There is only one such term, it
arises in the term $3D_{jk}U_{jk}(t_3)D_{jk}^{2}U_{jj}^\circ(t_2)$ of the
sum above, and by (\ref{ParU}) its contribution to $T_3^{(n)}$ is
\begin{equation*}
-\frac{\kappa _{4}i}{n}\sum_{k=1}^{n}\mathbf{E}\{ (U_{jj}\ast
U_{kk})(t_3)(U_{jj}\ast U_{jj}\ast U_{kk})(t_{2})\}.
\end{equation*}%
In view of (\ref{WUjj<}), (\ref{Ujjv}), and (\ref{v12}) we can replace here
all $U_{jj}$ and $U_{kk}$ with $v$ in the limit $n\rightarrow\infty$, so
that we have
\begin{align}
\lim_{n\rightarrow\infty}T_3^{(n)}(t_3,t_2)=-{\kappa _{4}i} (v\ast
v)(t_3)(v\ast v\ast v)(t_{2}).  \label{T3}
\end{align}
Summarizing (\ref{Covn}), (\ref{e3clt4}), (\ref{T1}) -- (\ref{T2}), and  (%
\ref{T3})  we obtain the equation with respect to $Cov$ of (\ref{Cov}):
\begin{align}
Cov(t_1,t_2)+w^{2}\int_{0}^{t_1}dt_3\int_{0}^{t_3}v(t_4)Cov(t_3-t_4,t_2)dt_4
=A(t_1,t_2),  \label{Coveq}
\end{align}%
where
\begin{equation}
A(t_1,t_2)=-2w^{2}\int_{0}^{t_1}\Phi (t_3,t_{2})dt_3+{\kappa _{4}}(v\ast
v\ast v)(t_{2}) \int_{0}^{t_1}(v\ast v)(t_3)dt_3  \label{A}
\end{equation}%
and $\Phi$ is given by (\ref{Phi}). To solve (\ref{Coveq}) we use the
generalized Fourier transform with respect to $t_1$ (see Proposition \ref%
{p:Four}). Note, that equation (\ref{Coveq}) is of the form (\ref{intrel}),
corresponding to $\delta =0$ in (\ref{supde}), thus we can use formulas (\ref%
{solut}) -- (\ref{TtSt}) to write its solution. Since
\begin{equation}
\widetilde{v}(z):=-i\int_{0}^{\infty }e^{-itz}v(t)dt=(2w^{2})^{-1}(\sqrt{%
z^{2}-4w^{2}}-z),  \label{f}
\end{equation}%
with the branch that is determined by the asymptotic $\sqrt{z^{2}-4w^{2}}%
=z+O(z^{-1}),\;z\rightarrow \infty$, and $(z+w^{2}\widetilde{v}(z))^{-1}=-%
\widetilde{v}(z)$, then we have for $T$ of (\ref{TtSt}):
\begin{equation}
T(t)=\frac{i}{2\pi }\int_{L}e^{itz}\frac{dz}{z+w^{2}\widetilde{v}(z)}=-v(t).
\label{T}
\end{equation}%
Hence, the unique differentiable solution of (\ref{Coveq}) is given by
\begin{equation}
Cov(t_1,t_2)=-2w^2\int_{0}^{t_1}v(t_1-t_3){\Phi} (t_3,t_{2})dt_3+{\kappa _{4}%
}\prod_{j=1}^2(v\ast v\ast v)(t_{j}).  \label{C12}
\end{equation}
We have by  (\ref{Phi}), (\ref{vt}), and a little algebra:
\begin{align*}
\int_{0}^{t_1}v(t_1-t_3){\Phi} (t_3,t_{2})dt_3=-\int_{[-2w,2w]^{3}}\frac{%
e^{i\lambda _{1}t_1}-e^{i\lambda _{2}t_1}}{\lambda _{1}-\lambda _{2}}\frac{%
e^{i\lambda _{2}t_2}-e^{i\lambda _{3}t_2}}{\lambda _{2}-\lambda _{3}}%
\prod_{j=1}^3\rho _{sc}(\lambda _{j})d\lambda _{j},
\end{align*}%
and by (\ref{f})
\begin{align*}
(v\ast v\ast v)(t)&=\frac{i}{2\pi }\int_{L}e^{itz}(2w^{2})^{-3}(\sqrt{%
z^{2}-4w^{2}}-z)^3dz \\
&=w^{-4}\int_{-2w}^{2w}e^{it\lambda}(w^2-\lambda^2) \rho
_{sc}(\lambda)d\lambda.
\end{align*}%
Hence, putting these expressions in (\ref{C12}), and then plugging the
result in (\ref{nCov}), we finally get:
\begin{align}
\lim_{n\rightarrow \infty }n\mathbf{Cov}\{&(\varphi_1 (M))_{jj}, (\varphi_2
(M))_{jj}\}  \notag \\
&=-2w^{2}\int_{[-2w,2w]^{3}}\frac{(\varphi_{1} (\lambda _{1})-\varphi_{1}
(\lambda _{2}))(\varphi_{2} (\lambda _{2})-\varphi_{2} (\lambda _{3}))}{%
(\lambda _{1}-\lambda _{2})(\lambda _{2}-\lambda _{3})}\prod_{j=1}^3\rho
_{sc}(\lambda _{j})d\lambda _{j}  \notag \\
&\quad+\kappa_4\prod_{j=1}^2w^{-4}\int_{-2w}^{2w}\varphi_j(\lambda)(w^2-%
\lambda^2) \rho _{sc}(\lambda)d\lambda .  \label{covph}
\end{align}%
Writing the numerator in the first integral as
\begin{equation}
\varphi_{1} (\lambda _{1})\varphi_{2} (\lambda _{2})-\varphi_{1} (\lambda
_{1})\varphi_{2} (\lambda _{3})-\varphi _{1}(\lambda _{2})\varphi_{2}
(\lambda _{2})+\varphi_{1} (\lambda _{2})\varphi_{2} (\lambda _{3}),
\label{phipro}
\end{equation}%
we observe that there is at least one integration which does not involve $%
\varphi $'s. This and the relation
\begin{equation*}
\int \frac{\rho _{sc}(\mu )d\mu }{\mu -\lambda }=-\lambda /2w^{2}
\end{equation*}%
allow us to deduce (\ref{Coff}) from (\ref{covph}). A simple way to perform
the corresponding calculations is to write the r.h.s. of (\ref{covph}) as
the limit as $\varepsilon \rightarrow 0$ of the same expression in which $%
\lambda _{3}$ is replaced by $\lambda _{3}+i\varepsilon $. One can also use
the Poincar\'{e} - Bertrand formula \cite{Mu:53} to deal with double
singular integrals, appearing after plugging (\ref{phipro}) in (\ref{covph}).
\end{proof}

\subsection{Limit Theorem for Matrix Elements}

\begin{theorem}
\label{t:clt} Consider the real symmetric Wigner random matrix of the form%
\begin{equation}
M_{n}=n^{-1/2}W_{n},\quad W_{n}=\{W_{jk}\in \mathbb{R},\;W_{jk}=W_{kj}=(1+%
\delta _{jk})^{1/2}V_{jk}\}_{j,k=1}^{n},  \label{MW1}
\end{equation}%
where $\{V_{jk}\}_{1\leq j\leq k<\infty }$  are i.i.d. random variables such
that
\begin{equation*}
\mathbf{E}\{V_{11}\}=0,\quad \mathbf{E}\{V_{11}^{2}\}=w^{2},\quad
\end{equation*}%
and the functions $\ln \mathbf{E}\{e^{itV_{11}}\}$ are entire.

Then for any $\varphi :\mathbb{R\rightarrow R}$ whose Fourier transform (\ref%
{FT}) satisfies (\ref{F3<}) and for any $j=j_{n}\in \lbrack 1,n]$ the random
variable $\sqrt{n}\varphi _{j_{n}j_{n}}^{\circ }(M)$ converges in
distribution as $n\rightarrow \infty $ to the random variable $\xi $ having
the characteristic function
\begin{equation}\label{ln}
\mathbf{E}\{e^{ix\xi }\}=\exp \{(-{x^{2}}V_{d}^{W}[\varphi ]+w^{2}x^{\ast
2})/2\}f(x^{\ast }),
\end{equation}%
where $f(x)=\mathbf{E}\{e^{ixV_{11}}\},$

\begin{equation}
x^*=\frac{\sqrt{2}x}{w^2}\int_{-2w}^{2w}\varphi(\mu)\mu\rho_{sc}(\mu)d\mu,
\label{x}
\end{equation}
$\rho_{sc}$ is the density of the semicircle law (\ref{rhosc}), and  $%
V_d^W[\varphi]$ is given by (\ref{VWd}).
\end{theorem}

\begin{remark}
Condition $W_{jk}=(1+\delta _{jk})^{1/2}V_{jk}$ is pure technical. In
particular, it can be shown that in the case of matrix $\widetilde{M}%
=n^{-1/2}V$, the Theorem \ref{t:clt} holds true with
\begin{equation*}
x^*=\frac{x}{w^2}\int_{-2w}^{2w}\varphi(\mu)\mu\rho_{sc}(\mu)d\mu,
\end{equation*}
and
\begin{equation*}
\mathbf{E}\{e^{ix\xi }\}=\exp \{(-{x^{2}}V_{d}^{W}[\varphi ]+2w^{2}x^{\ast
2})/2\}f(x^{\ast }).
\end{equation*}%
\end{remark}

\begin{proof}
Note first that in view of (\ref{beta}) and (\ref{MW1}) we can write
\begin{equation}
W_{jk}=\beta_{jk}^{-1/2}V_{jk}.  \label{Wjk}
\end{equation}
Besides, since  $\ln\mathbf{E}\{e^{itV_{11}}\}$ is entire then we have
\begin{align}
\sum_{l=1}^\infty\frac{x^l|\kappa_{l+1}|}{l!}<\infty,\quad \forall x>0, \quad
\label{ser<}
\end{align}
where $\kappa_l$ is the $l$th cumulant of $V_{11}$. We also have
\begin{equation}
\ln f(x) =-w^{2}x^{2}/2+\sum_{l=3}^\infty\frac{\kappa_{l}(ix)^l}{l!}
\label{lnf11}
\end{equation}%
(see (\ref{lt})).

We consider the characteristic functions
\begin{equation}
Z_{jn}(x)=\mathbf{E}\left\{ e^{ix\sqrt{n}\varphi _{jj}^{\circ }(M)}\right\}
\label{Znxj}
\end{equation}%
and prove that for any $x\in \mathbb{R}$
\begin{equation}
\lim_{n\rightarrow \infty }Z_{jn}(x)=\mathbf{E}\{e^{ix\xi}\}=:Z_{d}(x),
\label{chfGaj}
\end{equation}%
i.e., $Z_{d}(x)$ is given by the r.h.s. of (\ref{ln}).

Assume first that the Fourier transform (\ref{FT}) of $\varphi $ satisfies
\begin{equation}
\int |\widehat{\varphi }(t)||t|^{l}dt<C_{\varphi}l! \quad\forall l\in
\mathbb{N} ,  \label{phil<}
\end{equation}%
in particular, $\varphi $ is analytic in $|z|<1$. Since $Z_{jn}(0)=1$ and $%
Z_{jn}(x)$ is continuous, we can write the relation
\begin{equation}
Z_{jn}(x)=1+\int_{0}^{x}Z_{jn}^{\prime }(y)dy,\quad x\in \mathbb{R,}
\label{ZnZnpj}
\end{equation}%
showing that it suffices to prove that the sequence $\{Z_{jn}^{\prime }\}$
is uniformly bounded on any finite interval and that for any converging
subsequences $\{Z_{jn_{i}}\}_{i\geq 1}$ and $\{Z_{jn_{i}}^{\prime }\}_{i\geq
1}$ there exists $Z(x)$, such that
\begin{equation}
\lim_{i\rightarrow \infty }Z_{jn_{i}}(x)=Z(x),  \label{lz1}
\end{equation}
and
\begin{align}
\lim_{i\rightarrow \infty }Z_{jn_{i}}^{\prime }(x)=Z(x)\Big[%
-xV_d^W[\varphi]+\sum_{l=3}^\infty\frac{\kappa_{l}x^{l-1}}{(l-1)!} \Big(%
\frac{i\sqrt{2}}{w^2}\int_{-2w}^{2w}\varphi(\mu)\mu\rho_{sc}(\mu)d\mu\Big)^l%
\Big].  \label{lz2}
\end{align}
Indeed, if yes, then $Z(x)$ is a continuous function, satisfying for every $%
x\in\mathbb{R}$ the equation
\begin{align}
Z(x)=1-\int_{0}^{x}Z(y)\Big[&-yV_d^W[\varphi]  \notag \\
&+\sum_{l=3}^\infty\frac{\kappa_{l}y^{l-1}}{(l-1)!} \Big(\frac{i\sqrt{2}}{w^2%
}\int_{-2w}^{2w}\varphi(\mu)\mu\rho_{sc}(\mu)d\mu\Big)^l\Big]dy,
\label{eqzga}
\end{align}%
whose unique solution is the r.h.s. of (\ref{ln}).

We denote
\begin{equation}
e_{jn}(x)=e^{ix\sqrt{n}\varphi _{jj}^{\circ }(M)},\quad   \label{ejn}
\end{equation}%
and write according to (\ref{phijj}) and (\ref{Znxj})%
\begin{equation}
Z_{jn}^{\prime }(x)=i\mathbf{E}\left\{ \sqrt{n}\varphi _{jj}^{\circ }(M)e^{ix%
\sqrt{n}\varphi _{jj}^{\circ }(M)}\right\}= i\int \widehat{\varphi }%
(t)Y_{jn}(x,t)dt,  \label{dZY}
\end{equation}%
where
\begin{align}
&Y_{jn}(x,t)=\sqrt{n}\mathbf{E}\{U_{jj}(t)e^{\circ }_{jn}(x)\},  \label{Yjn}
\end{align}%
and $U$ is defined in (\ref{U}). It follows from the Schwarz inequality and (%
\ref{WUjj<}) that
\begin{equation}
|Y_{jn}(x,t)|\leq C(1+|t|)^3.  \label{Y<}
\end{equation}%
This and (\ref{phil<}) with $l=2$ yield that the sequence $Z^{\prime }_{jn}$
is uniformly bounded. Hence, there is a convergent subsequence $%
Z_{jn^{\prime }}$, and by the dominated convergence theorem to find its
limit as $n\rightarrow\infty$ it suffices to find the pointwise limit of the
corresponding  subsequence $Y_{jn^{\prime }}$.

Let us show now that sequences $\{\partial Y_{jn}/\partial x\}$ and $%
\{\partial Y_{jn}/\partial t\}$ are uniformly bounded in $(t,x)\in K\subset%
\mathbb{R}^2_+$, $n\in\mathbf{N}$, for any bounded $K$, so that the sequence
$\{Y_{jn}\}$ is equicontinuous on any finite set of $\mathbb{R}^2_+$, and
contains convergent subsequences.

Since $\overline{Y_{jn}(x,t)}=Y_{jn}(-x,-t)$, we can confine ourselves to
the half-plane $\mathbb{R}^2_+=\{t\geq 0,\,\,x\in \mathbb{R}\}$, and from
now on $t>0$.

It follows from (\ref{phijj}) that%
\begin{equation*}
\frac{\partial }{\partial x}Y_{jn}(x,t)=i\int \widehat{\varphi }(t_{1})n%
\mathbf{E}\{U_{jj}^{\circ }(t_1)U_{jj}^{\circ }(t)e_{jn}(x)\}dt_1,
\end{equation*}%
where by (\ref{WUjj<}) and the Schwarz inequality
\begin{equation*}
n|\mathbf{E}\{U_{jj}^{\circ }(t_1)U_{jj}^{\circ }(t)e_{jn}(x)\}|\leq n%
\mathbf{Var}^{1/2}\{U_{jj}(t_1)\}\mathbf{Var}^{1/2}\{U_{jj}(t)\}\leq
C(1+|t|)^3(1+|t_1|)^3.
\end{equation*}%
Hence, in view of (\ref{phil<}) the sequence $\{{\partial Y_{jn}}/{\partial x%
}\}$ is uniformly bounded.

We have also%
\begin{equation}
\frac{\partial }{\partial t}Y_{jn}(x,t)=i\sqrt{n}\mathbf{E}%
\{(MU)_{jj}(t)e^{\circ }_{jn}(x)\}=i\sum_{k=1}^n\mathbf{E}%
\{W_{jk}\Phi_{jk}(x,t)\},  \label{dtY}
\end{equation}%
where
\begin{equation}
\Phi_{jk}(x,t)=U_{jk}(t)e^{\circ}_{jn}(x).  \label{Phijk}
\end{equation}
To transform the r.h.s. of (\ref{dtY}) and show its boundedness, we apply an
analog of integration by parts formula proposed in Lemma \ref{l:difgen}.
Note that $D_{jk}^l\Phi_{jk}=O(n^{l/2})$ as $n\rightarrow\infty$, hence,
there is no such finite $p\in \mathbb{N}$ that $\varepsilon_p$ of (\ref%
{difgen}) vanishes as $n\rightarrow\infty$, and so we need infinite version
of "integration by parts formula" given by (\ref{difinf}). We will apply (%
\ref{difinf}) to every term of the r.h.s. of (\ref{dtY}), and to do this we
check first that $\Phi_{jk}(x,t)$ satisfies condition (\ref{al}). Indeed,
using the Leibnitz rule %
we obtain
\begin{align}
D_{jk}^l\Phi_{jk}(x,t)=\sum_{m=0}^{l} \Big(%
\begin{array}{ll}
l &  \\
m &
\end{array}%
\Big) D_{jk}^{l-m}U_{jk}(t)D_{jk}^me^{\circ}_{jn}(x),  \label{DjkF=}
\end{align}
where
\begin{equation}
D_{jk}^me_{jn}(x)=D_{jk}^{m-1}\big(ix\sqrt{n}e_{jn}(x)D_{jk}\varphi_{jj}(M)%
\big),  \label{De=}
\end{equation}
(see (\ref{ejn}) ), so that
\begin{align*}
&D_{jk}^me_{jn}(x)=e_{jn}(x)\sum_{r=1}^m(ix\sqrt{n})^r\sum_ {%
\begin{array}{ll}
q=(q_1,...,q_r): &  \\
q_1+...+q_r=m &
\end{array}%
} C_{q,r}\prod_{s=1}^rD_{jk}^{q_s}\varphi_{jj}(M),
\end{align*}
and
\begin{equation*}
\sum_{q,r}C_{q,r}\leq 2^m.
\end{equation*}
Hence,
\begin{equation*}
|D_{jk}^me_{jn}(x)|\leq\big(2\sqrt{n}(1+|x|)\big)^m \max_ {1\leq r\leq
m,\;\sum_{s=1}^r q_s=m} \prod_{s=1}^r|D_{jk}^{q_s}\varphi_{jj}(M)|,
\end{equation*}
where we have in view of (\ref{ocdlu}) 
and (\ref{phil<})
\begin{equation}
|D_{jk}^{q_s}\varphi_{jj}(M)|\leq\int|\widehat{\varphi}%
(\theta)||D_{jk}^{q_s}U_{jj}(\theta)|d\theta\leq C_\varphi2^{q_s},
\label{Dphi<}
\end{equation}
so that%
\begin{equation}
|D_{jk}^me_{jn}(x)|\leq\big(4C_\varphi\sqrt{n}(1+|x|)\big)^m.  \label{Dem<}
\end{equation}
This, (\ref{ocdlu}), and  (\ref{DjkF=}) yield
\begin{equation}
|D^l_{jk}\Phi_{jk}(x,t)|\leq (4C_\varphi\sqrt{n}(1+|x|+t))^l,\quad x\in
\mathbb{R},\; t>0.  \label{Phijk<}
\end{equation}
Thus, $\Phi_{jk}(x,t)$ satisfies condition (\ref{al}). Applying (\ref{difinf}%
) to every term of the r.h.s. of (\ref{dtY}) and taking into account (\ref%
{Wjk}), we get:
\begin{eqnarray}
\frac{\partial }{\partial t}Y_{jn}(x,t)=i\sum_{l=1}^\infty \frac{\kappa_{l+1}%
}{l!}S_l^{(n)},\quad S_l^{(n)}(x,t)=\frac{1}{(\sqrt{n})^l}\sum_{k=1}^n
\beta_{jk}^{-(l+1)/2} \mathbf{E}\{D_{jk}^l\Phi_{jk}(x,t)\}.  \label{dtY=}
\end{eqnarray}%
Let us show that this series converges uniformly in $(t,x)\in K$, $n\in%
\mathbb{N}$. In view of (\ref{ser<}) it suffices to show that
\begin{equation}
|S_l^{(n)}|\leq (C_K)^l, \quad \forall (t,x)\in K,\;n\in \mathbb{N},
\label{ckl}
\end{equation}
where $C_K$ is an absolute constant depending only on $K$.  Since
\begin{align}
S_l^{(n)}(x,t)=&\frac{2^{(l+1)/2}-1}{(\sqrt{n})^l}\mathbf{E}%
\{D_{jj}^l\Phi_{jj}(x,t)\}+\frac{1}{(\sqrt{n})^l}\sum_{k=1}^n \mathbf{E}%
\{D_{jk}^l\Phi_{jk}(x,t)\}  \label{sln}
\end{align}
(see (\ref{beta})), where in view of (\ref{Phijk<}) the first term of the
r.h.s. is bounded, it suffices to prove (\ref{ckl}) for the second term of
the r.h.s. of (\ref{sln}).

Using the Leibnitz rule, we write for $l\geq 2$
\begin{align}
\frac{1}{(\sqrt{n})^l}\sum_{k=1}^n D_{jk}^l\Phi_{jk}(x,t)=& \frac{1}{(\sqrt{n%
})^l}\sum_{k=1}^n U_{jk}D_{jk}^le_{jn}+\frac{l}{(\sqrt{n})^l}%
\sum_{k=1}^nD_{jk}U_{jk}D_{jk}^{l-1}e_{jn}  \label{Sl=} \\
&+\frac{1}{(\sqrt{n})^l}\sum_{k=1}^n\sum_{m=0}^{l-2} \Big(%
\begin{array}{ll}
l &  \\
m &
\end{array}%
\Big) %
D_{jk}^{l-m}U_{jk}D_{jk}^me_{jn}=:a_{l1}^{(n)}+a_{l2}^{(n)}+a_{l3}^{(n)},
\notag
\end{align}
where (cf (\ref{DjkF=}) -- (\ref{Phijk<}))

\begin{align}
|a_{l3}^{(n)}|\leq(4C_\varphi(1+|x|+t)\big)^{l}.  \label{al3<}
\end{align}
Applying (\ref{De=}) and then the Leibnitz rule again, we obtain for $%
a_{l1}^{(n)}$ of (\ref{Sl=}):
\begin{align}
a_{l1}^{(n)}=&\frac{ix}{(\sqrt{n})^{l-1}}\sum_{k=1}^n
U_{jk}(t)D_{jk}\varphi_{jj}(M)D_{jk}^{l-1}e_{jn}  \label{al1=} \\
&+\frac{ix}{(\sqrt{n})^{l-1}}\sum_{k=1}^n U_{jk}(t)\sum_{m=0}^{l-2} \Big(%
\begin{array}{ll}
l-1 &  \\
m &
\end{array}%
\Big) D_{jk}^{l-1-m}\varphi_{jj}(M)D_{jk}^me_{jn},  \notag
\end{align}
where the sum over $m$ is  bounded by
\begin{equation*}
(\sqrt{n})^{l-2}(4C_\varphi(1+|x|+t)\big)^{l-1}
\end{equation*}
(cf (\ref{Phijk<}) and (\ref{al3<})). Taking into account that $
|\sum_{k=1}^n U_{jk}(t)|\leq n^{-1/2}$, we see that the second
term of the r.h.s. of (\ref{al1=}) is bounded by
\begin{equation*}
(8C_\varphi(1+|x|+t)\big)^{l}.
\end{equation*}
Besides, it follows from (\ref{phijj}) and (\ref{ParU}) that
\begin{align}
&D_{jk}\varphi_{jj}(M)=2i\beta_{jk}\int\widehat{\varphi}(\theta)\int_0^%
\theta U_{jj}(\theta-\theta_1)U_{jk}(\theta_1)d\theta_1d\theta,
\label{Dphi=} \\
&|D_{jk}\varphi_{jj}(M)|\leq 2\int_0^\infty|\widehat{\varphi}%
(\theta)|\int_0^\theta |U_{jk}(\theta_1)|d\theta_1d\theta\leq 2C_\varphi.
\label{Df<}
\end{align}
This, (\ref{Dphi<}) -- (\ref{Dem<}), and (\ref{SUjk<}) allow us to show that
the first term of the r.h.s. of (\ref{al1=}) is bounded by
\begin{align*}
{2|x|\big(4C_\varphi(1+|x|)\big)^{l-1}} \int_0^\infty|\widehat{\varphi}%
(\theta)|\int_0^\theta
\sum_{k=1}^n|U_{jk}(t)|&|U_{jk}(\theta_1)|d\theta_1d\theta \\
&\leq(4C_\varphi(1+|x|+|t|)\big)^{l}.
\end{align*}
Hence,
\begin{equation}
|a_{l1}^{(n)}|\leq(4C_\varphi(1+|x|+t)\big)^{l}.  \label{al1<}
\end{equation}
Finally, since by (\ref{ParU}) -- (\ref{beta})
\begin{align*}
D_{jk}U_{jk}(t)=i
(U_{jj}*U_{kk}+U_{jk}*U_{jk})(t)-i\delta_{jk}(U_{jj}*U_{jj})(t),
\end{align*}
we have for $a_{l2}^{(n)}$ of (\ref{Sl=}):
\begin{align*}
a_{l2}^{(n)}=&\frac{il}{(\sqrt{n})^{l}}\sum_{k=1}^n
(U_{jj}*U_{kk})(t)D_{jk}^{l-1}e_{jn}+\frac{il}{(\sqrt{n})^{l}}\sum_{k=1}^n
(U_{jk}*U_{jk})(t)D_{jk}^{l-1}e_{jn} \\
&-\frac{il}{(\sqrt{n})^{l}} (U_{jj}*U_{jj})(t)D_{jj}^{l-1}e_{jn},
\end{align*}
where the last two terms are bounded by $l(1+t)(4C_\varphi(1+|x|)\big)^{l-1}$
in view of (\ref{Dem<}), (\ref{SUjk<}), and the bound $|(U_{ab}*U_{cd})(t)|%
\leq |t|$. Besides, it follows from (\ref{De=}) and (\ref{Df<}) that the
first term is bounded by:
\begin{align*}
\frac{lt}{(\sqrt{n})^{l}}\sum_{k=1}^n \Big|D_{jk}^{l-1}e_{jn}\Big|&=\frac{lt%
}{(\sqrt{n})^{l-1}}\sum_{k=1}^n \Big|D_{jk}\varphi_{jj}D_{jk}^{l-2}e_{jn}+
\sum_{m=0}^{l-3} \Big(%
\begin{array}{ll}
l-2 &  \\
m &
\end{array}%
\Big) D_{jk}^{l-2-m}\varphi_{jj}D_{jk}^{m}e_{jn}\Big| \\
&\leq \frac{lt|x|}{(\sqrt{n})^{l-1}} \bigg[2 \int_0^\infty|\widehat{\varphi}%
(\theta)|\int_0^\theta \sum_{k=1}^n|U_{jk}(\theta_1)|d\theta_1d\theta\big(%
4C_\varphi\sqrt{n}(1+|x|)\big)^{l-2}  \notag \\
&\quad+(\sqrt{n})^{l-1}(4C_\varphi(1+|x|+|t|)\big)^{(l-2)}\bigg]\leq
(4C_\varphi(1+|x|+t)\big)^{l} ,
\end{align*}
so that
\begin{equation}
|a_{l2}^{(n)}|\leq (4C_\varphi(1+|x|+t)\big)^{l}.  \label{al2<}
\end{equation}%
Now (\ref{ckl}) with $l\geq 2$ follows from (\ref{sln}) -- (\ref{al3<}) and (%
\ref{al1<}) -- (\ref{al2<}). Hence, the series in (\ref{dtY}) converges
uniformly in $(t,x)\in K\subset\mathbb{R}^2_+$ and $n\in\mathbb{N}$.

To prove the boundedness of the sequence $\{\partial Y_{jn}/\partial t\}$,
it remains to make sure that $S_{1}^{(n)}$ is bounded. Applying (\ref{ParU})
-- (\ref{beta}) and (\ref{De=}), we obtain
\begin{align}
S_{1}^{(n)} (x,t)&=n^{-1/2}\sum_{k=1}^n \beta_{jk}^{-1} \mathbf{E}%
\{D_{jk}\Phi_{jk}(x,t)\}  \notag \\
&=in^{-1/2}\mathbf{E}\big\{\big(n(U_{jj}*v_n)(t)+tU_{jj}(t)\big)%
e_{jn}^\circ(x)\big\}  \notag \\
&\quad-2x\int\widehat{\varphi}(\theta)\int_0^\theta \mathbf{E}%
\{U_{jj}(\theta-\theta_1)U_{jj}(t+\theta_1)e_{jn}(x)\}d\theta_1d\theta
\label{S1=}
\end{align}
where $v_n$ is defined in (\ref{vn=}). Writing
\begin{align}
n^{1/2}\mathbf{E}\{U_{jj}(t_1)&v_n(t_2)e_{jn}^\circ(x)\}  \notag \\
&=^{}\mathbf{E}\{v_n(t_2)\}Y_{jn}(x,t_1)+ n^{1/2}\mathbf{E}%
\{U_{jj}(t_1)v_n^\circ(t_2)e_{jn}^\circ(x)\},  \label{S1Uu=}
\end{align}
and taking into account bounds (\ref{Y<}), (\ref{Varu<}), $|U_{jj}|\leq 1$,
and $|e_{jn}^\circ(x)|\leq 2$, we conclude that the r.h.s. of (\ref{S1Uu=})
is bounded, and so does $S_1^{(n)}$. Hence, the sequence $\{\partial
Y_{jn}/\partial t\}$ is uniformly bounded in $(t,x)\in K\subset\mathbb{R}^2_+
$, $n\in\mathbb{N}$.

Now it follows from the above  that the sequence $\{Y_{jn}\}$ is
equicontinuous on any bounded set of $\mathbb{R}^{2}$. Hence, for any
converging subsequence $\{Z_{jn_i}\}$ (see (\ref{lz1})) there is a
converging subsequence $\{Y_{jn^{\prime }_i}\}$  and function $Y$ (which
obviously depends on $\{Z_{jn_i}\}$) such that
\begin{equation}
\lim_{n^{\prime }_i\rightarrow\infty}Y_{jn^{\prime }_i}=Y,\quad
\lim_{n^{\prime }_i\rightarrow\infty}Z_{jn^{\prime }_i}=Z.  \label{YZ}
\end{equation}
We will show now that $Y$ satisfies certain integral equation leading
through (\ref{dZY}) to (\ref{eqzga}), hence, to (\ref{ln}).  This will
finish the proof of the theorem under condition (\ref{phil<}).

Applying the Duhamel formula (\ref{Duh}) and then (\ref{dtY}) and (\ref{dtY=}%
), we obtain
\begin{align}
Y_{jn}(x,t) =i\sqrt{n}\int_{0}^{t}\sum_{k=1}^{n}\mathbf{E}%
\{M_{jk}\Phi_{jk}(x,t_1)\}dt_1 =i\int_0^t\sum_{l=1}^\infty \frac{\kappa_{l+1}%
}{l!}S_l^{(n)}(x,t_1)dt_1,  \label{Yjn=}
\end{align}%
where $\Phi_{jk}$ and $S_l^{(n)}$ are defined in (\ref{Phijk}) and (\ref%
{dtY=}), respectively. In view of the uniform convergence of the series, to
make the limiting transition as $n\rightarrow\infty$ it suffices to find the
limits
\begin{align*}
S_l=\lim_{n\rightarrow\infty}S_l^{(n)}
\end{align*}%
for every fixed $l\in\mathbb{R}$.

Let us start with $S_1^{(n)}$. It follows from (\ref{S1=}) -- (\ref{S1Uu=})
that (cf (\ref{T1=}) -- (\ref{r}))
\begin{align}
S_{1}^{(n)}(x,t_1)=&i\int_0^{t_1}\mathbf{E}\{v_n(t_1-t_2)\}Y_{jn}(x,t_2)dt_2
\label{S1r1} \\
&-2xZ_{jn}(x)\int\widehat{\varphi}(\theta)\int_0^\theta \mathbf{E}%
\{U_{jj}(\theta-\theta_1)\}\mathbf{E}\{U_{jj}(t_1+\theta_1)\}d\theta_1d%
\theta +r_n(x,t_1),  \notag
\end{align}
where
\begin{align*}
r_n(x,t)=&in^{-1/2}\mathbf{E}\big\{\big[n(U_{jj}*v^\circ_n)(t)+tU_{jj}(t)%
\big]e_{jn}^\circ(x)\big\} \\
&-2x\int\widehat{\varphi}(\theta)\int_0^\theta \mathbf{E}\big\{%
U_{jj}(\theta-\theta_1)\big[U_{jj}^{\circ}(t+\theta_1)Z_{jn}(x)+U_{jj}(t+%
\theta_1)e_{jn}^\circ(x)\big]\big\}d\theta_1d\theta,
\end{align*}
and by (\ref{VarUv}) and boundedness of $U_{ab}$ and $e_{jn}$ we have
\begin{equation}
r_{n}=O(n^{-1/2}),\quad n\rightarrow\infty.  \label{rn1}
\end{equation}%
%
This, (\ref{limU}), and (\ref{YZ}) imply
\begin{align}
\lim_{n^{\prime }_i\rightarrow\infty} S_{1}^{(n^{\prime
}_i)}=&i\int_0^{t_1}v(t_1-t_2)Y(x,t_1)dt_2 -2xZ(x)\int\widehat{\varphi}%
(\theta)\Phi(t_1,\theta)d\theta,  \label{limS1}
\end{align}
where $\Phi$ is given by (\ref{Phi}).

In general case we have for $S_{l}^{(n)}$, $l\geq 2,$ of (\ref{dtY=}):
\begin{align}
S_{l}^{(n)}=&\frac{1}{(\sqrt{n})^l}\sum_{k=1}^n \beta_{jk}^{-(l+1)/2}\mathbf{%
E}\Big\{U_{jk}D_{jk}^le_{jn}+lD_{jk}U_{jk}D_{jk}^{l-1}e_{jn}  \notag \\
&+{l(l-1)}D_{jk}^2U_{jk}D_{jk}^{l-2}e_{jn}/2+(1-\delta_{l2})%
\sum_{m=0}^{l-3} \Big(%
\begin{array}{ll}
l &  \\
m &
\end{array}%
\Big) D_{jk}^{l-m}U_{jk}D_{jk}^me_{jn}\Big\}  \notag \\
=&S_{l1}^{(n)}+S_{l2}^{(n)}+S_{l3}^{(n)}+(1-\delta_{l2})S_{l4}^{(n)},
\label{S1234}
\end{align}
and we have (cf (\ref{al3<}))
\begin{equation}
|S_{l4}^{(n)}|\leq (4C_\varphi(1+|x|+|t|)\big)^{l}n^{-1/2}, \quad l\geq 3.
\label{Sl4}
\end{equation}

Now we use the rule that in fact has been used several times before and
which follows from (\ref{SUjk<}) and the boundedness of $e_{jn}$ of (\ref%
{ejn}): the presence of a single factor $U_{ak}$ in terms of the sum $%
\sum_{k=1}^n$ is equivalent to the presence of the factor $n^{-1/2}$, and
the presence of two or more factors $U_{ak}$, $U_{bk}$,... is equivalent to
the presence of the factor $n^{-1}$. It follows from (\ref{ParU}) that all
terms of $D_{jk}^2U_{jk}$ contain $U_{jk}$, besides, we have
\begin{equation}
D_{jk}^{l-2}e_{jn}=O((\sqrt{n})^{l-2}),\quad n\rightarrow\infty  \label{el2}
\end{equation}
for $l>2$, hence,
\begin{equation}
S_{l3}^{(n)}=O(n^{-1/2}),\quad n\rightarrow\infty.  \label{Sl3}
\end{equation}
We also have (see (\ref{al1=})):
\begin{align*}
S_{l1}^{(n)}=\frac{ix}{(\sqrt{n})^{l-1}}\sum_{k=1}^n \beta_{jk}^{-(l+1)/2}%
\mathbf{E}\big\{&U_{jk}(t_1)\big[D_{jk}\varphi_{jj}(M)D_{jk}^{l-1}e_{jn} \\
&+(l-1)D_{jk}^2\varphi_{jj}(M)D_{jk}^{l-2}e_{jn}\big]\big\}+O(n^{-1/2}),\quad
n\rightarrow\infty,
\end{align*}
where $D_{jk}\varphi_{jj}(M)$ is given by (\ref{Dphi=}), and
\begin{align}
D_{jk}^2\varphi_{jj}(M)=-2i\beta^2_{jk}\int\widehat{\varphi}(\theta)
(3U_{jj}*U_{jk}*U_{jk}+U_{jj}*U_{jj}*U_{kk})(\theta)d\theta,  \label{D2phi=}
\end{align}
Hence, using again the above rule, we get
\begin{align}
S_{l1}^{(n)}=&-2x\int\widehat{\varphi}(\theta)\frac{1}{(\sqrt{n})^{l-1}}%
\sum_{k=1}^n \beta_{jk}^{-(l-1)/2}\mathbf{E}\big\{%
U_{jk}(t_1)(U_{jj}*U_{jk})(\theta)D_{jk}^{l-1}e_{jn}  \notag \\
&+i(l-1)U_{jk}(t_1)(U_{jj}*U_{jj}*U_{kk})(\theta)D_{jk}^{l-2}e_{jn}%
\big\}d\theta  \notag \\
&+O(n^{-1/2}),\quad n\rightarrow\infty.  \label{Sl1=}
\end{align}
It follows from (\ref{jkkk1}) -- (\ref{jkkk2}) that if $l=2$ then the
contribution of the second term of the r.h.s. of (\ref{Sl1=}) vanishes as $%
n\rightarrow\infty$. To show this in the case $l\geq 3$ we use the evident
relation
\begin{equation}
D_{jk}^{m}e_{jn}(x)=(ix\sqrt{n}D_{jk}\varphi_{jj}(M))^me_{jn}(x)+O((\sqrt{n}%
)^{m-1}),\quad n\rightarrow\infty,  \label{Dme}
\end{equation}
with $D_{jk}\varphi_{jj}(M)$ given by (\ref{Dphi=}), and conclude that the
second term contains either $U_{jk}U_{jk}\cdot O((\sqrt{n})^{l-2})$ or $%
U_{jk}\cdot O((\sqrt{n})^{l-3})$, hence, we have in view of (\ref{SUjk<}):
\begin{align*}
n^{-(l-1)/2}\sum_{k=1}^n
\beta_{jk}^2U_{jk}(t_1)(U_{jj}*U_{jj}*U_{kk})(\theta)D_{jk}^{l-2}e_{jn}(x)=
O(n^{-1/2}),\quad n\rightarrow\infty.
\end{align*}
The first term of the r.h.s. of (\ref{Sl1=}) already contains $U_{jk}U_{jk}$%
. Thus its non-vanishing is due to the term $(ix\sqrt{n}%
D_{jk}\varphi_{jj}(M))^{l-1}e_{jn}(x)$ of $D_{jk}^{l-1}e_{jn}(x)$, and we
get:
\begin{align*}
S_{l1}^{(n)}&=\sum_{k=1}^n\beta_{jk}^{-(l+1)/2}\mathbf{E}\big\{ %
U_{jk}(t_1)(ixD_{jk}\varphi_{jj}(M))^{l}e_{jn}(x)\big\}+O(n^{-1/2}) \\
&=\sum_{k=1}^n\mathbf{E}\Big\{e_{jn}(x)U_{jk}(t_1)\Big(-2x\int\widehat{%
\varphi}(\theta) (U_{jk}*U_{jj})(\theta)d\theta\Big)^l\Big\} \\
&\quad\quad+(2^{(1-l)/2}-1)\mathbf{E}\Big\{e_{jn}(x)U_{jj}(t_1) \Big(-2x\int%
\widehat{\varphi}(\theta) (U_{jj}*U_{jj})(\theta)d\theta\Big)^l\Big\} \\
&\quad\quad+O(n^{-1/2}),\quad n\rightarrow\infty,
\end{align*}
where we took into account (\ref{Dphi=}) and the equality $%
\beta_{jj}^{(l-1)/2}=2^{(1-l)/2}$ (see (\ref{beta})). Applying (\ref{vn1}), (%
\ref{Ujjv}), and (\ref{v2}) we get for $l\ge 2$:%
\begin{align*}
\mathbf{E}\Big\{e_{jn}(x)\sum_{k=1}^nU_{jk}(t_1)&\prod_{m=1}^lU_{jk}(%
\theta_m)U_{jj}(\tau_m)\Big\} \\
&=Z_{jn}(x) \overline{v}_{n2}(t_1,\theta_1,...,\theta_m)\prod_{m=1}^l\mathbf{E%
}\{U_{jj}(\tau_m)\}+O(n^{-1/4}) \\
&=Z(x)v(t_1)\prod_{m=1}^lv(\theta_m)v(\tau_m)+o(1),\quad
n_{i}\rightarrow\infty,
\end{align*}
and
\begin{equation*}
\mathbf{E}\Big\{e_{jn}(x)\prod_{m=1}^{l^{\prime }}U_{jj}(\tau_m)\Big\}=
Z(x)\prod_{m=1}^{l^{\prime }}v(\tau_m)+o(1),\quad n_{i}\rightarrow\infty.
\end{equation*}
Besides, we have in view of (\ref{vt})
\begin{equation}
(v*v)(\theta)=-\frac{i}{w^2}\int_{-2w}^{2w}e^{i\mu\theta}\mu\rho_{sc}(\mu)d%
\mu.  \label{vv}
\end{equation}
The above allows to write%
\begin{align}  \label{Sl1}
\lim_{i\rightarrow\infty}S_{l1}^{(n_i)}=\sqrt{2}Z(x)v(t_1)\Big(\frac{i\sqrt{2%
}x}{w^2}\int_{-2w}^{2w}\varphi(\mu)\mu\rho_{sc}(\mu)d\mu\Big)^l.
\end{align}
It remains to analyze $S_{l2}^{(n)}$ of (\ref{S1234}):
\begin{align*}
S_{l2}^{(n)}&=\frac{l}{(\sqrt{n})^l}\sum_{k=1}^n \beta_{jk}^{-(l+1)/2}%
\mathbf{E}\{D_{jk}U_{jk}D_{jk}^{l-1}e_{jn}\} \\
&=-\frac{lx}{(\sqrt{n})^{l-1}}\sum_{k=1}^n\beta_{jk}^{-(l-1)/2} \mathbf{E}%
\{(U_{jj}*U_{kk}+U_{jk}*U_{jk})(t_1)D_{jk}^{l-2}(e_{jn}(x)D_{jk}%
\varphi_{jj}(M))\} \\
&=-\frac{lx}{(\sqrt{n})^{l-1}}\sum_{k=1}^n \mathbf{E}\Big\{%
(U_{jj}*U_{kk})(t_1)\Big(D_{jk}^{l-2}e_{jn}(x)D_{jk}\varphi_{jj}(M) \\
&\quad\quad\quad\quad\quad\quad\quad\quad\quad+(1-%
\delta_{l2})(l-2)D_{jk}^{l-3}e_{jn}(x)D_{jk}^2\varphi_{jj}(M)\Big)\Big\} %
+O(n^{-1/2}),\quad n\rightarrow\infty,
\end{align*}
where we used consequently (\ref{ParU}), (\ref{De=}), and then the Leibnitz
rule, (\ref{Dem<}), and (\ref{SUjk<}). Treating the first term of the last
expression analogously to the second term of (\ref{Sl1=}), we see that it is
of the order $O(n^{-1/2})$, $n\rightarrow\infty$. Hence, taking into account
(\ref{D2phi=}) and (\ref{SUjk<}), we get $S_{22}=0,$ and for $l\geq 3$
\begin{align*}
S_{l2}^{(n)}=2l(l-2)x&\int\widehat{\varphi}(\theta)\frac{1}{(\sqrt{n})^{l-1}}%
\sum_{k=1}^n \mathbf{E}\Big\{(U_{jj}*U_{kk})(t_1) \\
&\times(U_{jj}*U_{jj}*U_{kk})(\theta)D_{jk}^{l-3}e_{jn}(x)\Big\}%
d\theta+O(n^{-1/2}),\quad n\rightarrow\infty.
\end{align*}
If $l>3$, then in view of (\ref{Dme}), $D_{jk}^{l-3}e_{jn}(x)$ gives either
factor $U_{jk}O((\sqrt{n})^{l-3})$ or $O((\sqrt{n})^{l-4})$, that in both
cases leads to
\begin{equation}
S_{l2}^{(n)}=O(n^{-1/2}),\quad n\rightarrow\infty, \quad l>3.  \label{S42}
\end{equation}
If $l=3$, then
\begin{align*}
S_{32}^{(n)}=&6x\int\widehat{\varphi}(\theta)\frac{1}{n}\sum_{k=1}^n \mathbf{%
E}\{(U_{jj}*U_{kk})(t_1) (U_{jj}*U_{jj}*U_{kk})(\theta)e_{jn}(x)\}d\theta \\
&+O(n^{-1/2}),\quad n\rightarrow\infty,
\end{align*}
and it follows from (\ref{vn12}), (\ref{Ujjv}), and (\ref{v12}) that%
\begin{align*}
\mathbf{E}\{e_{jn}U_{jj}(\tau_1)&U_{jj}(\tau_2)U_{jj}(\tau_3)n^{-1}%
\sum_{k=1}^nU_{kk}(\tau_4)U_{kk}(\tau_5)\} \\
&=Z_{jn}(x)\prod_{m=1}^3\mathbf{E}\{U_{jj}(\tau_m)\}\mathbf{E}%
\{n^{-1}\sum_{k=1}^nU_{kk}(\tau_4)U_{kk}(\tau_5)\} +O(n^{-1/2}) \\
&=Z(z)\prod_{m=1}^5v(\tau_m)+o(1), \quad n_i\rightarrow\infty.
\end{align*}
Hence,
\begin{align*}
\lim_{n_i\rightarrow\infty}S_{32}^{(n_i)}=6xZ(x)(v*v)(t_1)\int\widehat{%
\varphi}(\theta)(v*v*v)(\theta)d\theta.
\end{align*}
This, (\ref{S1234}) -- (\ref{Sl3}), and (\ref{Sl1}) -- (\ref{S42}) yield for
$l\geq 2$:
\begin{align}
\lim_{n_i\rightarrow\infty}S_{l}^{(n_i)}=&Z(x)v(t_1)\Big(\frac{i\sqrt{2}x}{%
w^2}\int_{-2w}^{2w} \varphi(\mu)\mu\rho_{sc}(\mu)d\mu\Big)^l  \notag \\
&+\delta_{3l}6xZ(x)(v*v)(t_1)\int\widehat{\varphi}(\theta)(v*v*v)(\theta)d%
\theta.  \label{limS3}
\end{align}
Summarizing (\ref{Yjn=}), (\ref{limS1}), and (\ref{limS3}) we see that $Y$
of (\ref{YZ}) satisfies the equation
\begin{align}
Y(x,t)+w^{2}\int_{0}^{t}dt_1&\int_{0}^{t_1}v(t_1-t_2) Y(x,t_2)dt_2
\label{eqydxt} \\
& =ixZ(x)\int\widehat{\varphi}(\theta)A(t,\theta)d\theta+i\sqrt{2}%
B(x)Z(x)\int_{0}^{t}v(t_1)dt_1,  \notag
\end{align}%
where $A$ is defined in (\ref{A}), and
\begin{equation}
B(x)=\sum_{l=2}^\infty\frac{\kappa_{l+1}}{l!} \Big(\frac{i\sqrt{2}x}{w^2}%
\int_{-2w}^{2w}\varphi(\mu)\mu\rho_{sc}(\mu)d\mu\Big)^l.  \label{B}
\end{equation}
The kernel of this equation coincides with that of (\ref{Coveq}). Hence, the
argument leading to (\ref{C12}) and based on using of generalized Fourier
transform, yields
\begin{align*}
Y(x,t)=ixZ(x) \int\widehat{\varphi}(\theta)Cov(t,\theta)d\theta+i\sqrt{2}%
B(x)Z(x)(v*v)(t),
\end{align*}%
where $Cov$ is given by (\ref{C12}). Plugging this expression in (\ref{dZY})
and taking into account (\ref{vv}) and the equality (see (\ref{nCov}))%
\begin{align*}
\int\int\widehat{\varphi}(\theta)\widehat{\varphi}(t)Cov(t,\theta)dtd%
\theta=V_d^W[\varphi],
\end{align*}%
we finally get for $Z$ of (\ref{YZ}):
\begin{align*}
Z^{\prime }(x)=Z(x)\left[-xV_d^W[\varphi]+B(x)\frac{i\sqrt{2}}{w^2}%
\int_{-2w}^{2w}\varphi(\mu)\mu\rho_{sc}(\mu)d\mu\right].
\end{align*}%
This yields (\ref{lz1}) -- (\ref{lz2}), thus proves the theorem under
condition (\ref{phil<}).

The case of $\varphi\in E=\{\psi:\int(1+|t|)^{3}|\widehat{\psi }%
(t)|dt<\infty\}$ can be obtained via a standard approximation procedure.
Indeed, since the set $D=\{\varphi:\int|\widehat{\varphi\ }%
(t)||t|^{l}dt<C_\varphi l!,\;\forall l\in\mathbb{N}\}$ is big enough (in
particular, it contains functions $e^{-x^2}P_m(x)$, where $P_m(x)$ is a
polynomial), then for any  $\varphi \in E$ there exists a sequence $%
\{\varphi _{k}\}\subset D$, such  that
\begin{equation}
\lim_{k\rightarrow \infty }\int_{-2w}^{2w}|\varphi (\lambda )-\varphi
_{k}(\lambda )|d\lambda=0.  \label{fktof}
\end{equation}%
Denote for the moment the characteristic functions of (\ref{Znxj}) and (\ref%
{chfGaj}) as $Z_{n}[\varphi ]$ and $Z[\varphi ]$, to make explicit their
dependence on  $\varphi $. We have then for any $\varphi \in E$
\begin{align}
|Z_{n}[\varphi ]-Z[\varphi ]|&\leq |Z_{n}[\varphi ]-Z_{n}[\varphi
_{k}]|+|Z_{n}[\varphi _{k}]-Z[\varphi _{k}]|+|Z[\varphi _{k}]-Z[\varphi ]|
\notag \\
&:=T_{nk}^{(1)}+T_{nk}^{(2)}+T_{nk}^{(3)}.  \label{chaZ}
\end{align}%
The second term of the r.h.s. vanishes after the limit $n\rightarrow \infty $
in view of the above proof, since $\varphi _{k}\in D$. For the first term we
have from (\ref{Znxj}) and the Schwarz inequality that
\begin{eqnarray*}
|T_{nk}^{(1)}| &\leq &|x|\big(n\mathbf{Var}^{}\{(\psi_k(M))_{jj}\}\big)%
^{1/2},\quad \psi_k=\varphi-\varphi_k,
\end{eqnarray*}%
and then Theorem \ref{t:Cov} implies that
\begin{eqnarray*}
\limsup_{n\rightarrow \infty }|T_{nk}^{(1)}| &\leq &|x|(V_d^W[\psi_k])^{1/2}.
\end{eqnarray*}%
Since $V^W_{d}$ of (\ref{VWd}) is continuous with respect to the $L^{1}$
convergence, then in view of (\ref{fktof}) $T_{nk}^{(1)}$ vanishes after the
subsequent limits $n\rightarrow \infty$, $k\rightarrow \infty $.

At last, we have by (\ref{ser<}) and the continuity of the r.h.s.  of (\ref%
{ln}) with respect to the $L^{1}$ convergence, that the third term of (\ref%
{chaZ}) vanishes after the limit $k\rightarrow \infty .$ Thus, we have
proved the Central Limit Theorem under condition (\ref{F3<}).
\end{proof}

\begin{corollary}
It follows from Theorem \ref{t:clt} that if $\varphi$ is even, then the
random variable $\sqrt{n}\varphi ^{\circ }(M)_{jj}$ converges in
distribution to the Gaussian random variable with zero mean and the variance
$V_d^W[\varphi]$.
\end{corollary}

\begin{remark}
Random matrix theory deals mostly with eigenvalues of large random matrices. However, the statistical properties of eigenvectors are also of considerable interests for a number of reasons, in particular in view of possible links with the problem of existence of absolutely continuous spectrum of the multidimensional Schr\'odinger operator with random potential (see e.g. \cite{Er:10r,Kh-Pa:93}). In the case of the Gaussian random matrices (GOE, null Wishart) the eigenvectors are rotationally invariant and according to recent works \cite{Ba-Co:07,Er:10r,Le-Pe:09} the eigenvectors of the non-Gaussian random matrices (Wigner, sample covariance) are similar in several aspects to the eigenvectors of the Gaussian random matrices. On the other hand, the results of \cite{L-Pa:09} and this papers imply that there are asymptotic properties of eigenvectors of the non-Gaussian random matrices which are different of those for  the Gaussian random matrices.
\end{remark}

\subsection{Auxiliary results}

\begin{lemma}
\label{l:main} Consider the unitary matrix $U(t)=e^{itM}$ of (\ref{U}) -- (%
\ref{norU}), where $M$ is the Wigner matrix (\ref{MW}) -- (\ref{Wmom12}),
define
\begin{align}
&v_n(t)=n^{-1}\sum_{k=1}^nU_{kk}(t),  \notag \\
&v_n(t_1,t_2)=n^{-1}\sum_{k=1}^nU_{kk}(t_1)U_{kk}(t_2),  \label{vn12} \\
&v_{n1}(t_1,t_2)=n^{-1/2}\sum_{k=1}^nU_{jk}(t_1) U_{kk}(t_2),\quad
\label{vn2} \\
&v_{n2}(\overline{t})=\sum_{k=1}^n\prod_{m=1}^l U_{jk}(t_m),\quad\; l\geq 3,
\label{vn1}
\end{align}
where $\overline{t}=(t_1,..,t_l)$, and put $\overline{f}=\mathbf{E}\{f\}$.
Then we have:
\begin{equation}
\lim_{n\rightarrow\infty}\overline{v}_n(t)=\int
_{-2w}^{2w}e^{it\lambda}\rho_{sc}(\lambda)d\lambda=:v(t),  \label{vt}
\end{equation}
where $\rho _{sc}$ is the density of the semicircle law (\ref{rhosc}), and
under the conditions of Theorem \ref{t:Cov}
\begin{alignat}{3}
&\text{(i)} & &\mathbf{Var}\{U_{jj}(t)\} \leq C(1+|t|)^6/n, & &
\lim_{n\rightarrow\infty}\overline{U}_{jj}(t)=v(t),  \label{Ujjv} \\
&\text{(ii)} & &\mathbf{Var}\{v_n(t)\}=O(n^{-2}),\;n\rightarrow\infty, & &
\label{Varu<} \\
&\text{(iii)} & &\mathbf{Var}\{v_n(t_1,t_2)\}=O(n^{-1/2}),\;n\rightarrow%
\infty, & & \lim_{n\rightarrow\infty}\overline{v}_n(t_1,t_2)=v(t_1)v(t_2),
\label{v12} \\
& \text{(iv)} & &\mathbf{Var}\{v_{n1}(t_1,t_2)\}=O(n^{-1/2}),\;n\rightarrow%
\infty, & & \lim_{n\rightarrow\infty}\overline{v}_{n1}(t_1,t_2)=0,
\label{v1} \\
&\text{(v)} & &\mathbf{Var}\{v_{n2}(\overline{t})\}=O(n^{-1/2}),\;n%
\rightarrow\infty, & & \lim_{n\rightarrow\infty}\overline{v}_{n2}(\overline{t%
})=\prod_{m=1}^l v(t_m),  \label{v2}
\end{alignat}
where $O(n^\alpha)$ can depend on $t$, $\overline{t}$.
\end{lemma}

\begin{remark}
It can be shown that all statements of the lemma remain valid under
conditions of Theorem \ref{t:clt}
\end{remark}

\begin{proof}
Statement (\ref{vt}) follows from the well known fact of random matrix
theory  (see e.g. \cite{Pa:05} and references therein) according to which
for any bounded and continuous $\varphi$
\begin{equation*}
\lim_{n\rightarrow\infty}n^{-1}\mathbf{E}\{\Tr \varphi(M)\}=\int_{-2w}^{2w}
\varphi(\lambda)\rho_{sc}(\lambda)d\lambda,
\end{equation*}
where $\rho_{sc}$ is the density of the semicircle law (\ref{rhosc}).

(i) Let $\widehat{M}=n^{-1/2}\widehat{W}$ be GOE matrix (\ref{GOE})
independent of $M$, and

\begin{equation}
\widehat{U}(t)=e^{it\widehat{M}}.  \label{UGOE}
\end{equation}
We can write
\begin{align}
V_{n}&:=\mathbf{Var}\{U_{jj}(t)\} = \mathbf{E}\{(U_{jj}(t)-\widehat{U}%
_{jj}(t))U_{jj}^{\circ }(-t)\},  \label{K12}
\end{align}%
and then follow the interpolation procedure proposed in \cite{L-Pa:09}.
Namely, consider the "interpolating" random matrix (see \cite{Pa:07,L-Pa:09}%
)
\begin{equation}
M(s)=s^{1/2}M+(1-s)^{1/2}\widehat{M},\quad 0\leq s\leq 1,  \label{Mt}
\end{equation}%
viewed as defined on the product of the probability spaces of matrices $W$
and $\widehat{W}$. We denote again by $\mathbf{E}\{\dots \}$ the
corresponding expectation in the product space. Since $M(1)=M$, $M(0)=%
\widehat{M}$, then putting
\begin{equation}
U(t,s)=e^{itM(s)},  \label{Uts}
\end{equation}%
we obtain
\begin{align}
U_{jj}(t)-\widehat{U}_{jj}(t)&=\int_{0}^{1}\frac{\partial }{\partial s}%
U_{jj}(t,s) ds  \label{U-U} \\
&=\frac{i}{2}\int_{0}^{1}\sum_{l,m=1}^n \Big(\frac{1}{\sqrt{sn}}W^{(n)}_{lm}-%
\frac{1}{\sqrt{(1-s)n}}\widehat{W}_{lm}\Big)(U_{jl}*U_{mj})(t,s)ds.  \notag
\end{align}%
Thus,
\begin{align*}
V_n=\frac{i}{2}\int_{0}^{1}\bigg[&\frac{1}{\sqrt{sn}}\sum_{l,m=1}^n \mathbf{E%
}\{{W}^{(n)}_{lm}\Phi_{lm}\} -\frac{1}{\sqrt{(1-s)n}}\sum_{l,m=1}^n \mathbf{E%
}\{\widehat{W}_{lm}\Phi_{lm}\}\bigg]ds,
\end{align*}%
where

\begin{equation}
\Phi_{lm}=(U_{jl}*U_{mj})(t,s)U_{jj}^{\circ}(-t).  \label{Philm}
\end{equation}
A simple algebra based on (\ref{ParU}) and (\ref{ocdlu}) allows to obtain
\begin{equation}
|D_{lm}^q\Phi_{lm}|\leq C _{q}|t|^{q+1},  \label{P<}
\end{equation}
with $C_{q}$ depending only on $q\in\mathbb{N}$.

Now, applying differentiation formula (\ref{difgen}) with $p=4$ and $%
\Phi=\Phi_{lm}$ to every term of the first sum and differentiation formula (%
\ref{diffga}) to every term of the second sum, we obtain (cf (\ref{Covn}) --
(\ref{e3clt4})):
\begin{align}
V_{n}=\frac{i}{2}\int_{0}^{1}\Big[\sum_{l=2}^4 T_p^{(n)}+\varepsilon_4\Big]%
s^{-1/2}ds,  \label{K1}
\end{align}
where
\begin{equation}
T_p^{(n)}=\frac{1}{p!n^{(p+1)/2}}\sum_{l,m=1}^{n}\kappa _{p+1,lm}\mathbf{E}%
\big\{D_{lm}^{p}\Phi_{lm}\big\},\quad p=2,3,4,  \label{Tp}
\end{equation}
and by (\ref{b3}) and (\ref{P<})
\begin{equation}
|\varepsilon_4|\leq\frac{C_4w_6}{n^{3}}\sum_{l,m=1}^n\sup_{M\in \mathcal{S}%
_{n}}|D_{lm}^{5}\Phi_{lm}|\leq {C(1+|t|)^6}{n^{-1}}.  \label{e2<}
\end{equation}
Consider now $T^{(n)}_{2}$ and note that by (\ref{cums}) and (\ref{mu3}) $%
\kappa_{3,lm}=\mu_3$, so that
\begin{align}
T^{(n)}_{2}&=\frac{\mu_3}{2n^{3/2}}\sum_{l,m=1}^n \mathbf{E}%
\{D^2_{lm}\Phi_{lm}\}  \label{lT2=} \\
&=\frac{\mu_3}{2n^{3/2}}\sum_{l,m=1}^n \mathbf{E}\{U_{jj}^{%
\circ}(-t)D^2_{lm}(U_{jl}*U_{mj})(t,s)  \notag \\
&\hspace{2cm}+(U_{jl}*U_{mj})(t,s)D_{lm}^{2}U_{jj}(-t)  \notag \\
&\hspace{2cm}+2D_{lm}(U_{jl}*U_{mj})(t,s)D_{lm}U_{jj}(-t)\}:=%
\mu_3[T_{21}^{(n)}+T_{22}^{(n)}+T_{23}^{(n)}].  \notag
\end{align}
It follows from (\ref{ParU}) that
\begin{align}
T_{21}^{(n)}=-3\mu_3{s^{}n}^{-3/2}\sum_{l,m=1}^n \mathbf{E}%
\{&U_{jj}^{\circ}(-t)(U_{jl}*U_{jm}*U_{lm}*U_{lm}  \label{T21=} \\
&+2U_{jl}*U_{jl}*U_{lm}*U_{mm}+U_{jl}*U_{jm}*U_{ll}*U_{mm})(t,s)\}.  \notag
\end{align}
Here by (\ref{SUjk<})
\begin{align}
&\sum_{l,m=1}^n |U_{jl}U_{jm}U_{lm}U_{lm}|\leq 1,\quad \sum_{l,m=1}^n |
U_{jl}U_{jl}U_{lm}U_{mm}|\leq n^{1/2},  \label{good} \\
& \sum_{l,m=1}^n | U_{jl}U_{jm}U_{ll}U_{mm}|\leq n^{}.  \label{bad}
\end{align}
Hence, applying the Schwarz inequality and taking into account that $%
(1*1*1*1)(t)=t^3/6,$ we obtain
\begin{equation}
|T_{21}^{(n)}|\leq C(1+|t|^{3})(n^{-1/2}V^{1/2}_{n}+n^{-1}).  \label{T21<}
\end{equation}
The terms $T^{(n)}_{22}$ and $T_{23}^{(n)}$ contain sums of two types
\begin{align}
\sum_{l,m=1}^n U_{jl}U_{jl}U_{jl}U_{jm}U_{mm}=O(n^{1/2}),\quad
\sum_{l,m=1}^n U_{jl}U_{jl}U_{jm}U_{jm}U_{lm} =O(1),\quad n\rightarrow\infty,
\label{T23}
\end{align}
where the r.h.s. of both equalities follows from (\ref{SUjk<}). Hence, $%
|T_{22}^{(n)}+T_{23}^{(n)}|\leq C(1+|t|^{3})n^{-1}$. This, (\ref{lT2=}), and
(\ref{T21<}) yield
\begin{equation}
|T_{2}^{(n)}|\leq C(1+|t|^{3})(n^{-1/2}V^{1/2}_{n}+n^{-1}).  \label{T2<}
\end{equation}
Acting in the similar way and taking into account (\ref{kapmu}) and (\ref%
{w6<}) implying $|\kappa_{p,jk}|\leq C$, $p=4,5$, we get analogous bounds
for $T^{(n)}_{3}$ and $T^{(n)}_{4}$ of (\ref{K1}):
\begin{equation}
|T_{l}^{(n)}|\leq {C(1+|t|)^{l+1}}{n^{-1}},\quad l=4,5.  \label{Tl<}
\end{equation}
Note, that in the case of $T^{(n)}_{3}$ and $T^{(n)}_{4}$ the argument  is
even simpler because here we have factors $n^{-2}$ and $n^{-5/2}$,
respectively, instead of $n^{-3/2}$ of $T^{(n)}_{2}$.

Now it follows from (\ref{K1}) -- (\ref{e2<}) and (\ref{T2<}) -- (\ref{Tl<})
that $V_{n}^{1/2}$ satisfies the quadratic inequality:
\begin{equation*}
V_{n}-{C(1+|t|)^3}{n^{-1/2}}V_{n}^{1/2}-{C(1+|t|)^6}{n^{-1}}\leq 0,
\end{equation*}
implying (\ref{WUjj<}) and then (\ref{VarF}).

To finish the proof of (i) it remains to show that
\begin{equation}  \label{Uv}
\lim_{n\rightarrow\infty}\mathbf{E}\{U_{jj}(t)\}=v(t).
\end{equation}
In the GOE case we have by the orthogonal invariance of GOE probability
measure and (\ref{vt}):
\begin{equation}  \label{vGOE}
\mathbf{E}\{\widehat{U}_{jj}(t)\}=\mathbf{E}\{n^{-1}\Tr \widehat{U}%
(t)\}\rightarrow v(t) \;\text{as}\;n\rightarrow\infty.
\end{equation}
Besides, it follows from (\ref{Mt}), (\ref{U-U}), (\ref{diffga}), and (\ref%
{difgen}) with $p=3$, that
\begin{align*}
\mathbf{E}\{U_{jj}(t)-\widehat{U}_{jj}(t)\}=&\frac{i}{2}\int_{0}^{1}\Big[%
\sum_{l=2}^3 T_p^{(n)}+\varepsilon_3\Big]s^{-1/2}ds,
\end{align*}
where $T_p$ are given by (\ref{Tp}) with $\Phi_{lm}=(U_{jl}*U_{mj})(t,s)$
(cf (\ref{Philm})), and
\begin{equation*}
|\varepsilon_3|\leq\frac{C_3w^{5/6}_6}{n^{5/2}}\sum_{l,m=1}^n\sup_{M\in
\mathcal{S}_{n}}|D_{lm}^{4}\Phi_{lm}|\leq {C(1+|t|)^5}{n^{-1/2}}.
\end{equation*}
A similar but much simpler argument leading to (\ref{T2<}) -- (\ref{Tl<})
allows to conclude that $|T_p^{(n)}|\leq O(n^{-1/2})$, $n\rightarrow\infty$,
$p=2,3$. Hence,
\begin{align*}
\mathbf{E}\{U_{jj}(t)-\widehat{U}_{jj}(t)\}=O(n^{-1/2}),\;n\rightarrow\infty.
\end{align*}
This and (\ref{vGOE}) yield (\ref{Uv}) and finish the proof of (i).

(ii) The proof of (\ref{Varu<}) repeats with natural modifications the one
of the first part of (i). Namely, similarly to (\ref{K12}) -- (\ref{K1}) we
have for $V_{n}=\mathbf{Var}\{v_{n}(t)\}$:
\begin{align}
V_{n}&=\mathbf{E} \{[v_{n}(t)-\widehat{v}_{n}(t)]v^\circ_{n}(-t)\}  \notag \\
&=\frac{it}{2}\int_{0}^{1}\bigg[\frac{1}{\sqrt{sn^{3}}}\sum_{i,k=1}^n
\mathbf{E}\{{W}^{(n)}_{ik}\Phi_{ik}\} -\frac{1}{\sqrt{(1-s)n^{3}}}%
\sum_{i,k=1}^n \mathbf{E}\{\widehat{W}_{ik}\Phi_{ik}\}\bigg]ds  \notag \\
&=\frac{it}{2}\int_{0}^{1}\Big[\sum_{p=2}^4T^{(n)}_{p}+\varepsilon_4\Big]%
s^{-1/2}ds,  \label{Ett}
\end{align}
where now
\begin{align*}
&T^{(n)}_{p}=\frac{1}{p!n^{(3+p)/2}}\sum_{i,k=1}^n \kappa _{p+1,ik}\mathbf{E}%
\{D^p_{ik}\Phi_{ik}\},\quad p=2,3,4,  \notag \\
&\Phi_{ik}=U_{k}(t,s){v}^\circ_{n}(-t),\quad |D_{ik}^l(s)\Phi_{ik}|\leq
C(|t|),
\end{align*}
$C(|t|)$ is a polynomial in $|t|$ with positive coefficients, and
\begin{align}
&|\varepsilon _{4}|\leq \frac{C_{4}w^{}_6}{n^{4}}\sum_{i,k=1}^{n}\sup_{M\in
\mathcal{S} _{n}}\Big|D_{ik}^{5}\Phi_{ik}\Big|\leq C(|t|) n^{-2}.
\label{e3<}
\end{align}
Using the argument leading to (\ref{T2<}) and (\ref{Tl<}) it can be shown
that
\begin{equation*}
T_p^{(n)}\leq C(|t|)(n^{-1}V_{n}^{1/2}+n^{-2}),\;p=2,3, \quad
T_{4}^{(n)}\leq C(|t|)n^{-2}.
\end{equation*}
This and (\ref{Ett}) -- (\ref{e3<}) allow us to write the inequality
\begin{equation*}
V_{n}\leq C(|t|)(n^{-1}V_{n}^{1/2}+n^{-2}),
\end{equation*}
valid for any real $t$ and implying (\ref{Varu<}).

(iii) Statement (iii) was proved in Lemma 3.1 of \cite{Ly-Pa:08}.

(iv) Let $\overline{t}=(t_1,t_2)$ and
\begin{equation*}
\widehat{v}_{n1}(\overline{t})=n^{-1/2}\sum_{k=1}^n \widehat{U}_{jk}(t_1)%
\widehat{U}_{kk}(t_2)
\end{equation*}
with $\widehat{U}$ of (\ref{UGOE}). We have similar to (\ref{U-U})
\begin{align}
v_{n1}(\overline{t})-\widehat{v}_{n1}(\overline{t})&=\int_{0}^{1}n^{-1/2}%
\frac{\partial }{\partial s}\sum_{k=1}^nU_{jk}(t_1,s)U_{kk}(t_2,s)ds  \notag
\\
&=\frac{i}{2}\int_{0}^{1}n^{-1/2}\sum_{p,q=1}^n \big(%
s^{-1/2}W^{(n)}_{pq}-(1-s)^{-1/2}\widehat{W}_{pq}\big)\Phi_{pq}(\overline{t}%
,s)ds,  \label{v1v1}
\end{align}
where now%
\begin{align*}
&\Phi_{pq}(\overline{t},s)=n^{-1/2}%
\sum_{k=1}^n(U_{jp}*U_{qk})(t_1,s)U_{kk}(t_2,s)+U_{jk}(t_1,s)(U_{kp}*U_{qk})(t_2,s),
\\
&|D_{pq}^{i}\Phi_{pq}(\overline{t},s)|\leq C_{}(\overline{t}), \quad i\in%
\mathbb{N}.
\end{align*}
Applying (\ref{diffga}) and (\ref{difgen}) with $p=3$ and $\Phi=\Phi_{pq}$,
we get
\begin{align}
\mathbf{E}\{v_{n1}(\overline{t})-\widehat{v}_{n1}(\overline{t}%
)\}=\int_{0}^{1} \bigg[&\frac{\mu_3}{2n^{3/2}}\sum_{p,q=1}^n \mathbf{E}\big\{%
D_{pq}^{2}\Phi_{pq}(\overline{t},s)\big\}  \label{v-v} \\
&+\frac{\kappa_{4}}{6n^{2}}\sum_{p,q=1}^n \mathbf{E}\big\{%
D_{pq}^{3}\Phi_{pq}(\overline{t},s)\big\}\bigg] s^{-1/2}ds+
O(n^{-1/2}),\quad n\rightarrow\infty,  \notag
\end{align}
where $\kappa_{4}$ is defined in (\ref{k4}). Using the argument leading to (%
\ref{T2<}) and (\ref{Tl<}) and based on (\ref{norU}) -- (\ref{ParU}), it can
be shown that both the terms in the square brackets of (\ref{v-v}) are of
the order $O(n^{-1/2})$, $n\rightarrow\infty$.  Hence,
\begin{align}
\mathbf{E}\{v_{n1}(\overline{t})\}&-\mathbf{E}\{\widehat{v}_{n1}(\overline{t}%
)\} =O(n^{-1/2}),\quad n\rightarrow\infty.  \label{vn1vn1}
\end{align}
Moreover, replacing $\Phi_{pq}(\overline{t},s)$ with $\Phi_{pq}(\overline{t}%
,s){v}_{n1}^{\circ}(\overline{t})$ we can also obtain
\begin{align}
\mathbf{Var}\{v_{n1}(\overline{t})\}=\mathbf{E}\{(v_{n1}(\overline{t})-%
\widehat{v}_{n1}(\overline{t})){v}_{n1}^{\circ}(\overline{t})\}
=O(n^{-1/2}),\quad n\rightarrow\infty.  \label{varvn1}
\end{align}
Now it follows from (\ref{vn1vn1}) -- (\ref{varvn1}) that to finish the
proof of (iii) it remains to show that
\begin{equation}
\lim_{n\rightarrow\infty}\mathbf{E}\{\widehat{v}_{n1}(\overline{t})\}=0.
\label{limv2G}
\end{equation}
Indeed, since by the orthogonal invariance of the GOE probability measure we
have
\begin{equation*}
\mathbf{E}\{\widehat{U}_{jk}(t)\}=\delta_{jk}\mathbf{E}\{\widehat{v}%
_{n}(t)\},\quad
\end{equation*}
where $\widehat{v}_{n}(t)=\Tr
\widehat{U}(t)$, $|\widehat{v}_{n}(t)|\leq 1$, then in view of (\ref{WUjj<})
and (\ref{SUjk<})
\begin{align}
\mathbf{E}\{\widehat{v}_{n1}(t_1,t_2)\} &=n^{-1/2}\mathbf{E}\{\widehat{v}%
_{n}(t_{1})\}\mathbf{E}\{\widehat{v}_{n}(t_{2})\}  \label{vn1G} \\
&+n^{-1/2}\sum_{k=1}^n\mathbf{E}\Big\{\widehat{U}_{jk}(t_1) \widehat{U}%
^\circ_{kk}(t_2)\Big\}=O(n^{-1/2}),\quad n\rightarrow\infty,  \notag
\end{align}
and we get (iv).

(v) The scheme of the proof of (iv) is the same as the one of (iii). Namely,
we have similar to (\ref{v1v1}) -- (\ref{v-v}) (see also (\ref{U-U})):
\begin{align*}
v_{n2}(\overline{t})-\widehat{v}_{n2}(\overline{t}) =\frac{i}{2}%
\sum_{m=1}^l\int_{0}^{1}\frac{1}{\sqrt{n}}\sum_{p,q=1}^n \big(%
s^{-1/2}W^{(n)}_{pq}-(1-s)^{-1/2}\widehat{W}_{pq}\big)\Phi_{pqm}(\overline{t}%
,s)ds,
\end{align*}
where now
\begin{align*}
&\Phi_{pqm}(\overline{t},s)=\sum_{k=1}^n(U_{jp}*U_{qk})(t_m,s)\prod_{m^{%
\prime }\neq m}U_{jk}(t_{m^{\prime }},s), \\
&|D_{pq}^{l}\Phi_{pqm}(\overline{t},s)|\leq C_{}(\overline{t}), \quad l\in%
\mathbb{N}.
\end{align*}
Applying (\ref{diffga}) and (\ref{difgen}) with $p=3$ one can get an analog
of (\ref{v-v}) and then show that (cf (\ref{vn1vn1}) -- (\ref{varvn1}))
\begin{equation*}
\mathbf{E}\{v_{n2}(\overline{t})\}-\mathbf{E}\{\widehat{v}_{n2}(\overline{t}%
)\}= O(n^{-1/2}),\quad n\rightarrow\infty,
\end{equation*}
and
\begin{equation*}
\mathbf{Var}\{v_{n2}(\overline{t})\}=O(n^{-1/2}),\quad \quad
n\rightarrow\infty.
\end{equation*}
So, it remains to prove that
\begin{equation}
\lim_{n\rightarrow\infty}\mathbf{E}\{\widehat{v}_{n1}(\overline{t}%
)\}=\prod_{m=1}^l v(t_m).  \label{limv}
\end{equation}
Applying (\ref{diffga}) and (\ref{ParU}) and then (\ref{vn=}) and (\ref%
{SUjk<}), we get
\begin{align}
\frac{\partial}{\partial t_1}\mathbf{E}\{\widehat{v}_{n2}(\overline{t})\}&=%
\frac{i}{\sqrt{n}} \sum_{k,p=1}^n\mathbf{E}\Big\{\widehat{W}_{jp}\widehat{U}%
_{kp}(t_1) \prod_{m=2}^l\widehat{U}_{jk}(t_m) \Big\}   \label{dv=} \\
&=-\frac{w^2}{{n}} \sum_{k=1}^n\mathbf{E}\Big\{(t_1\widehat{U}_{jk}(t_1)+(n%
\widehat{v}_{n}*\widehat{U}_{jk})(t_1)) \prod_{m=2}^l\widehat{U}%
_{jk}(t_m)+\sum_{m=2}^l\prod_{m^{\prime }\neq m}\widehat{U}%
_{jk}(t_{m^{\prime }})  \notag \\
&\quad\quad\quad\quad\quad\quad\times \int_0^{t_m} \Big(\widehat{U}%
_{jj}(t_m-s_1)\widehat{U}_{kk}(t_1+s_1)+\widehat{U}_{jk}(t_m-s_1)\widehat{U}%
_{jk}(t_1+s_1)\Big)ds_1\Big\},  \notag
\end{align}
so that by (\ref{vn=}) and (\ref{SUjk<}) ${\partial \mathbf{E}\{\widehat{v}%
_{n2}(\overline{t})\}}/{\partial t_1}=O(1)$, $n\rightarrow\infty$, and by
the symmetry
\begin{equation*}
\frac{\partial}{\partial t_m}\mathbf{E}\{\widehat{v}_{n2}(\overline{t}%
)\}=O(1),\quad \quad n\rightarrow\infty, \quad m=1,..,l.
\end{equation*}
Hence, there exists a subsequence $\{\mathbf{E}\{{\widehat{v}}_{n_i2}(%
\overline{t})\}\}$ that converges uniformly on any compact set of $\mathbb{R}%
^l$. Now, applying the Duhamel formula (\ref{Duh}) and then (\ref{dv=}), we
obtain
\begin{align*}
\mathbf{E}\{\widehat{v}_{n2}(\overline{t})\}&=\mathbf{E}\Big\{\prod_{m=2}^l%
\widehat{U}_{jj}(t_m)\Big\} +\int_0^{t_1}\frac{i}{\sqrt{n}} \sum_{k,p=l}^n%
\mathbf{E}\Big\{\widehat{W}_{jp}\widehat{U}_{kp}(s) \prod_{m=2}^l\widehat{U}%
_{jk}(t_m) \Big\}ds \\
&=\mathbf{E}\Big\{\prod_{m=2}^l\widehat{U}_{jj}(t_m)\Big\}-{w^2}%
\int_0^{t_1}ds\int_0^s \mathbf{E}\{\widehat{v}_{n}(s-s_1)\}\mathbf{E}\{%
\widehat{v}_{n2}(s_1,t_2,...,t_l)\}ds_1+r_n,
\end{align*}
where
\begin{align*}
r_n=&-\frac{w^2}{{n}}\int_0^{t_1} \sum_{k=l}^n\mathbf{E}\Big\{(s\widehat{U}%
_{jk}(s)+(n\widehat{v}^\circ_{n}*\widehat{U}_{jk})(s)) \prod_{m=2}^l\widehat{%
U}_{jk}(t_m) \\
& +\sum_{m=2}^l\prod_{m^{\prime }\neq m}\widehat{U}_{jk}(t_{m^{\prime
}})\int_0^{t_m} (\widehat{U}_{jj}(t_m-s_1)\widehat{U}_{kk}(s+s_1)+\widehat{U}%
_{jk}(t_m-s_1)\widehat{U}_{jk}(s+s_1))ds_1\Big\}ds,
\end{align*}
and by the Schwarz inequality, (\ref{SUjk<}), and (\ref{Varu<})
\begin{equation*}
r_n=O(n^{-1}),\quad n\rightarrow\infty.
\end{equation*}
Taking into account (\ref{Ujjv}) we obtain that every limit of converging
subsequence
\begin{equation*}
v_2(\overline{t})=\lim_{i\rightarrow\infty}\widehat{v}_{n_i2}(\overline{t})
\end{equation*}
satisfies the equation:
\begin{align}
v_2(\overline{t})+w^{2}\int_{0}^{t_1}ds%
\int_{0}^{s}v(s-s_1)v_2(s,t_2,...,t_l)ds_1=\prod_{m=2}^lv(t_{m}).
\label{v1eq}
\end{align}
Applying the generalized Fourier transform with respect to the variable $t_1$
(see Proposition \ref{p:Four}), we get
\begin{equation*}
\widetilde{v_2}(z,t_2,...,t_l)(1+w^2z^{-1}\widetilde{v}(z))=z^{-1}%
\prod_{m=2}^lv(t_{m})
\end{equation*}
with $\widetilde{v}$ of (\ref{f}). Hence, $\widetilde{v_2}(z,t_2,...,t_l)=%
\widetilde{v}(z)\prod_{m=2}^lv(t_{m})$, and
\begin{equation*}
{v_2}(\overline{t})=\prod_{m=1}^lv(t_{m}).
\end{equation*}
This completes the proof of (v) and the proof of the Lemma.
\end{proof}

\bigskip

After this paper was completed we became aware of paper \cite{So-Co:11} by Soshnikov et al
in which Theorem \ref{t:clt} was proved by another method and under weaker conditions.


\begin{thebibliography}{99}


\bibitem{Ba-Co:07}
 Bai, Z. D., Miao, B. Q., Pan, G. M.: {On asymptotics of eigenvectors of
 large sample covariance matrix}. Ann. Probab. 35 (2007)  1532--1572.







\bibitem{Di-Co:03} D'Aristotile, A., Diaconis, P., Newman, C.: \emph{ Brownian
motion and the classical groups}. Probability, statistics and their
applications: papers in honor of Rabi Bhattacharya, IMS Lecture Notes
Monogr. Ser. 41, Inst. Math. Statist., Beachwood, OH, 97--116 (2003)


\bibitem{Er:10r}
Erd\"os, L.: {Universality of Wigner random matrices: a Survey of Recent Results}. arXiv:1004.0861




\bibitem{KKP:96} Khorunzhy, A.M., Khoruzhenko, B.A., Pastur, L.A.: {%
Asymptotic properties of large random matrices with independent entries}. J.
Math. Phys. {37} (1996) 5033--5060.

\bibitem{Kh-Pa:93}
Khorunzhy, A. M., Pastur, L. A.: {Limits of infinite interaction radius,
dimensionality and the number of components for random operators with
off-diagonal randomness}.  Comm. Math. Phys.  153  (1993) 605--646.

\bibitem{Le-Pe:09}
Ledoit, O., Pe\'che\', S.: {Eigenvectors of some large sample covariance matrix
ensembles}. arXiv:0911.3010.


\bibitem{Ly-Pa:08} Lytova, A., Pastur, L.: {Central Limit Theorem for linear
eigenvalue statistics of random matrices with independent entries}, Annals of
Probability  37 (2009) 1778--1840.


\bibitem{L-Pa:09} Lytova, A., Pastur, L.: {Fluctuations of matrix elements of
regular functions of Gaussian random matrices}. J. Stat. Phys. 134
(2009) 147--159.

\bibitem{Me:91} Mehta, L.: {Random Matrices}. Academic Press, New York (1991).


\bibitem{Mu:53} Muskhelishvili, N.I.: {Singular Integral Equations}.
Noordhoff, Groningen (1953)


\bibitem{So-Co:11}
O'Rourke, S., Renfrew, D.,  Soshnikov, A.: {On fluctuations of matrix entries
of regular functions of Wigner matrices with Non-Identically Distributed Entries},
arXiv:1104.1663

\bibitem{Pa:05} Pastur, L.: {A simple approach to the global regime of
Gaussian Ensembles of random matrices}. Ukrainian Math. J. {57} (2005) 936--966.


\bibitem{Pa:07} Pastur, L.: Eigenvalue Distribution of Random
Matrices. In \textit{Random Media 2000, J. Wehr (Ed.)} 95--206.
Wydawnictwa ICM, Warsaw (2007).



\bibitem{Pr-Ro:69} Prokhorov, Yu. V.,  Rozanov, Yu. A.: {Probability Theory%
}, Springer, Berlin (1969).

\bibitem{Ti:86} Titchmarsh, E.C.: {Introduction to the theory of Fourier
integrals.} Chelsea Publishing Co., New York (1986)
\end{thebibliography}
\end{document}